# RELATIONS BETWEEN INVASION PERCOLATION AND CRITICAL PERCOLATION IN TWO DIMENSIONS

By Michael Damron[1], Artëm Sapozhnikov[2] and
Bálint Vágvölgyi[3]

*Courant Institute, EURANDOM and Vrije Universiteit Amsterdam*

We study invasion percolation in two dimensions. We compare connectivity properties of the origin's invaded region to those of (a) the critical percolation cluster of the origin and (b) the incipient infinite cluster. To exhibit similarities, we show that for any $k \geq 1$, the $k$-point function of the first so-called pond has the same asymptotic behavior as the probability that $k$ points are in the critical cluster of the origin. More prominent, though, are the differences. We show that there are infinitely many ponds that contain many large disjoint $p_c$-open clusters. Further, for $k > 1$, we compute the exact decay rate of the distribution of the radius of the $k$th pond and see that it differs from that of the radius of the critical cluster of the origin. We finish by showing that the invasion percolation measure and the incipient infinite cluster measure are mutually singular.

**1. Introduction.** Self-organized criticality has become a subject of great interest in recent years. Although there is no general definition for it, we can say that a system or model has this property if the definition of the model requires no parameter, yet some characteristics of the model resemble those at criticality of a parametric model with a phase transition. One such model is *invasion percolation*, a stochastic growth model that mirrors

Received June 2008; revised February 2009.

[1]Supported by NSF Grant OISE-0730136 (Percolative and Disordered Systems: A U.S.-Brazil-Netherlands Based International Collaboration).

[2]Supported by Netherlands Organisation for Scientific Research (NWO) Grant 613.000.429.

[3]Supported by Netherlands Organisation for Scientific Research (NWO) Grant 639.033.201.

aspects of the critical Bernoulli percolation picture without tuning any parameter. The invasion model was introduced independently by two groups ([2] and [11]), who studied it numerically. The first mathematically rigorous study of invasion percolation appeared in [4]. Connections between the invasion cluster and critical Bernoulli percolation have been established in, for instance, [4, 6, 18, 21] and [22], using both heuristics and rigorous arguments. These results indicated so many parallels between the invaded region and the incipient infinite cluster that a question naturally arose: to what extent are these objects similar? This question was studied on the regular tree in [1]. It was shown that, although the invaded region and the incipient infinite cluster are locally similar, globally, they differ significantly. In this paper, we prove local similarities between critical Bernoulli clusters and certain invaded clusters (the ponds) in the plane. We also show that, globally, the invaded region and the incipient infinite cluster are essentially different.

In the remainder of this section, we define the invasion percolation model and, using results of [4], we introduce the ponds of the invasion. We then review results concerning relations between invasion percolation and critical Bernoulli percolation. Finally, we state the main results of the paper.

### 1.1. *The model.*

For simplicity, *we restrict ourselves here to the square lattice.* Invasion percolation can be similarly defined on other two-dimensional lattices and the results of this paper still hold for lattices which are invariant under reflection in one of the coordinate axes and under rotation about the origin by some angle in $(0, \pi)$. In particular, this includes the triangular and honeycomb lattices.

Although our results concern invasion in the plane, we give the definition of invasion percolation for $\mathbb{Z}^d$. Consider the hypercubic lattice $\mathbb{Z}^d$ with its set of nearest neighbor bonds $\mathbb{E}^d$. We denote edges by their endpoints, that is, we write $e = \langle x, y \rangle$ if the two endpoints of $e$ are $x$ and $y$. Letting $G = (V, E)$ be an arbitrary subgraph of $(\mathbb{Z}^d, \mathbb{E}^d)$, we define the outer edge boundary $\Delta G$ of $G$ as follows:

$$\Delta G = \{e = \langle x, y \rangle \in \mathbb{E}^d : e \notin E(G), \text{but } x \in G \text{ or } y \in G\}.$$

The first step is to assign independent random variables, uniformly distributed in $[0, 1]$, to each bond $e \in \mathbb{E}^d$. We denote these variables by $\tau_e$. Using them, we recursively define an increasing sequence $G_0, G_1, G_2, \ldots$ of connected subgraphs of the lattice. $G_0$ only contains the origin, with no edges. Once $G_i = (V_i, E_i)$ is defined, we select the edge $e_{i+1}$ that minimizes $\tau$ on $\Delta G_i$. We take $E_{i+1} = E_i \cup \{e_{i+1}\}$ and let $G_{i+1}$ be the graph induced by the edge set $E_{i+1}$. The graph $G_i$ is called the *invaded region at time $i$*, and the graph $\mathcal{S} = \bigcup_{i=0}^{\infty} G_i$ is called the *invasion percolation cluster* (IPC). Let $E_\infty = \bigcup_{i=0}^{\infty} E_i$.

Since we would like to compare Bernoulli percolation to the invasion, we use a well-known analogous definition of Bernoulli percolation that makes the coupling of the two models immediate. For any $p \in [0, 1]$, we say that an edge $e \in \mathbb{E}^d$ is $p$-*open* if $\tau_e < p$. It is obvious that the resulting random graph of $p$-open edges has the same distribution as the one obtained by declaring each edge of $\mathbb{E}^d$ open with probability $p$ and closed with probability $1 - p$, independently of the states of all other edges. The percolation probability $\theta(p)$ is the probability that the origin is in the infinite cluster of $p$-open edges. There is a critical probability $p_c = \inf\{p : \theta(p) > 0\} \in (0, 1)$. For general background on Bernoulli percolation, we refer the reader to [5].

It was shown in [4] that for all $p > p_c$, the invasion intersects the infinite $p$-open cluster with probability 1. In the case $d = 2$, this result immediately follows from the Russo–Seymour–Welsh theorem (see Section 11.7 in [5]). Furthermore, the definition of the invasion mechanism implies that if the invasion reaches the $p$-open infinite cluster for some $p$, then it will never leave this cluster. Combining these facts yields that if $e_i$ is the edge added at time $i$, then $\limsup_{i \to \infty} \tau_{e_i} = p_c$. From now on, we consider only $d = 2$. In this case, it is well known that $\theta(p_c) = 0$, which implies that for every $t > 0$, there is an edge $e(t)$ such that $e(t)$ is invaded after step $t$ and $\tau_{e(t)} > p_c$. The last two results give that $\hat{\tau}_1 = \max\{\tau_e : e \in E_\infty\}$ exists and is greater than $p_c$. Let $\hat{e}_1$ denote the edge at which the maximum value of $\tau$ is taken and assume that $\hat{e}_1$ is invaded at step $i_1 + 1$. Following the terminology of [15], we call the graph $G_{i_1}$ the *first pond* of the invasion and denote it $\hat{V}_1$. The edge $\hat{e}_1$ is called the *first outlet*. The second pond of the invasion is defined similarly. Note that the same argument as above implies that $\hat{\tau}_2 = \max\{\tau_{e_i} : e_i \in E_\infty, i > i_1\}$ exists and is greater than $p_c$. If we assume that $\hat{\tau}_2$ is taken on the edge $\hat{e}_2$ at step $i_2 + 1$, we call the graph $G_{i_2} \backslash G_{i_1}$ the *second pond* of the invasion and denote it $\hat{V}_2$. The further ponds $\hat{V}_k$ can be defined analogously.

The following interpretation gives a natural meaning to the ponds. Consider an infinite piece of land divided into square parcels. These parcels are separated by dikes whose heights are given by the values of independent random variables, uniformly distributed on $[0, 1]$. One of the parcels, called the parcel of the origin, contains an infinite source of water. First, the water level in the parcel of the origin rises until it reaches the height of the lowest adjacent dike and then spills over into the parcel on the other side of this dike. Next, the water level rises in both parcels until it reaches the height of the lowest dike on the boundary of the union of the two parcels, at which time a new parcel floods. The process continues indefinitely and, as time approaches infinity, an infinite region of land will flood. Consider the dual lattice of $\mathbb{Z}^2$, each dual edge having the $\tau$ value of its corresponding edge in the original lattice, identifying the dual edges with the dikes and the origin with the source of water. Each vertex of $\mathbb{Z}^2$ corresponds to exactly one parcel of land. It is evident from the invasion mechanism and from the way the flood spreads on the land that a parcel is flooded if and only if the corresponding vertex of $\mathbb{Z}^2$ is invaded. We now explain the meaning of the first pond in the flood setting. At step $i_1$, when the first outlet is invaded, the minimal $\tau$ value on the boundary of $G_{i_1}$ is that of $\hat{e}_1$. However, this is the edge with the largest $\tau$ value ever added to the invasion. This means that the invasion will never return to $G_{i_1}$, that is, no edge on $\Delta G_{i_1}$, other than $\hat{e}_1$, will be invaded. Therefore, after some time, all water will flow over the dike corresponding to $\hat{e}_1$ and the water level in each parcel of the first pond will be constant and equal to $\hat{\tau}_1$. The same argument shows that after some time, the water level in the second pond will become, and remain, $\hat{\tau}_2$, and so on.

Now that our model is defined, we review a few results that established connections between the invasion and the critical percolation models. To the best of our knowledge, the first paper with mathematically rigorous results in this area was [4], where it was shown, among other things, that the empirical distribution of the $\tau$ value of the invaded edges converges to the uniform distribution on $[0, p_c]$. Results on the fractal nature of the invaded region were also obtained in [4]. The authors showed that the region has zero volume fraction, given that there is no percolation at criticality, and that it has boundary-to-volume ratio $(1 - p_c)/p_c$. This corresponds to the asymptotic boundary-to-volume ratio for large critical clusters (see [10] and [14]). The above results indicate that a large proportion of the edges in the IPC belong to big $p_c$-open clusters.

An object that turns out to be closely related to the invaded region is the *incipient infinite cluster* (IIC). Loosely speaking, one can say that the IIC is the "infinite open cluster at criticality." The IIC can be constructed by conditioning on the origin being connected to a site at distance $n$ from the origin in critical percolation and by considering the cluster of the origin. If we let $n \to \infty$, an infinite cluster is obtained and this cluster is called the incipient infinite cluster. (Later in this paper, we will give the precise definition. For detailed results on the IIC, we refer the reader to [8].) Let $S_n$ be the number of invaded sites within a distance of at most $n$ from the origin. The scaling of the moments of $S_n$ as $n$ goes to infinity was obtained in [6] and [22], and it turned out to coincide with the scaling of the corresponding moments for the IIC. Another similarity established in [6] is concerned with the invasion picture far away from the origin: the invasion measure was shown to be locally the same as the IIC measure.

The diameter and volume of the first pond of the invasion were studied in [18, 19]. It was shown that the decay rates of their distributions coincide, respectively, with the decay rates of the distributions of the diameter and the volume of the critical cluster of the origin in Bernoulli percolation.

To the best of knowledge, the only paper to date concerned with the differences between the invasion model and critical percolation is [1]. The authors consider invasion percolation on regular trees. The scaling behavior of the $r$-point function and the volume of the invaded region at and below a given height can be explicitly computed. It is found that while the power laws of the scaling are the same for the invaded region and for the incipient infinite cluster, the scaling functions differ and, consequently, the two clusters behave differently. In fact, their laws are found to be mutually singular. Even though the arguments of [1] do not work for invasion in the plane, their results give a strong indication that, in spite of the presence of many similarities, the two objects are indeed different.

In this paper, we compare connectivity properties of the origin's invaded region to those of the critical percolation cluster of the origin and the IIC. In Theorems 1.1 and 1.2, we give the asymptotic behavior for the $k$-point function of the first pond. We continue to study the relation between the IPC and large $p_c$-open clusters in Theorems 1.3 and 1.4. We show that, for any $K$ and $N$, there are infinitely many ponds that contain at least $K$ disjoint $p_c$-open clusters of size at least $N$. We also show that, provided the radius of the first pond is larger than $N$, the first pond contains at least $K$ disjoint $p_c$-open clusters of size at least $N$ with probability bounded from below by a positive constant independent of $N$. For $k > 1$, we compute the exact decay rate of the distribution of the radius of the $k$th pond in Theorem 1.5. Unlike the decay rate of the distribution of the radius of the first pond [18], it is strictly different from that of the radius of the critical cluster of the origin. Finally, in Theorem 1.8, we show that the IPC measure and the IIC measure are mutually singular.

### 1.2. *Notation.*

In this section, we set out most of the notation and definitions used in the paper.

For $a \in \mathbb{R}$, we write $|a|$ for the absolute value of $a$ and, for a site $x = (x_1, x_2) \in \mathbb{Z}^2$, we write $|x|$ for $\max(|x_1|, |x_2|)$. For $n > 0$ and $x \in \mathbb{Z}^2$, let $B(x, n) = \{y \in \mathbb{Z}^2 : |y - x| \le n\}$ and $\partial B(x, n) = \{y \in \mathbb{Z}^2 : |y - x| = n\}$. We write $B(n)$ for $B(0, n)$ and $\partial B(n)$ for $\partial B(0, n)$. For $m < n$ and $x \in \mathbb{Z}^2$, we define the annulus $\mathrm{Ann}(x; m, n) = B(x, n) \setminus B(x, m)$. We write $\mathrm{Ann}(m, n)$ for $\mathrm{Ann}(0; m, n)$.

We consider the square lattice $(\mathbb{Z}^2, \mathbb{E}^2)$, where $\mathbb{E}^2 = \{(x, y) \in \mathbb{Z}^2 \times \mathbb{Z}^2 : |x - y| = 1\}$. Let $(\mathbb{Z}^2)^* = (1/2, 1/2) + \mathbb{Z}^2$ and $(\mathbb{E}^2)^* = (1/2, 1/2) + \mathbb{E}^2$ be the vertices and the edges of the dual lattice. For $x \in \mathbb{Z}^2$, we write $x^*$ for $x + (1/2, 1/2)$. For an edge $e \in \mathbb{E}^2$, we denote its ends, left (resp., right) or bottom (resp., top), by $e_x, e_y \in \mathbb{Z}^2$. The edge $e^* = (e_x + (1/2, 1/2), e_y - (1/2, 1/2))$ is called the *dual edge* to $e$. Its ends, bottom (resp., top) or left (resp., right), are denoted by $e_x^*$ and $e_y^*$. Note that, in general, $e_x^*$ and $e_y^*$ are not the same as $(e_x)^*$ and $(e_y)^*$. For a subset $\mathcal{K} \subset \mathbb{Z}^2$, let $\mathcal{K}^* = (1/2, 1/2) + \mathcal{K}$. We say that an edge $e \in \mathbb{E}^2$ is in $\mathcal{K} \subset \mathbb{Z}^2$ if both of its ends are in $\mathcal{K}$.

Let $(\tau_e)_{e \in \mathbb{E}^2}$ be independent random variables, uniformly distributed on $[0, 1]$, indexed by edges. We call $\tau_e$ the *weight* of an edge $e$. We define the weight of an edge $e^*$ as $\tau_{e^*} = \tau_e$. We denote the underlying probability measure by $\mathbb{P}$ and the space of configurations by $([0, 1]^{\mathbb{E}^2}, \mathcal{F})$, where $\mathcal{F}$ is a natural $\sigma$-field on $[0, 1]^{\mathbb{E}^2}$. We say that an edge $e$ is $p$-*open* if $\tau_e < p$ and $p$-*closed* if $\tau_e > p$. An edge $e^*$ is $p$-*open* if $e$ is $p$-open and it is $p$-*closed* if $e$ is $p$-closed. The event that two sets of sites $\mathcal{K}_1, \mathcal{K}_2 \subset \mathbb{Z}^2$ are connected by a $p$-open path is denoted by $\mathcal{K}_1 \overset{p}{\longleftrightarrow} \mathcal{K}_2$ and the event that two sets of sites $\mathcal{K}_1^*, \mathcal{K}_2^* \subset (\mathbb{Z}^2)^*$ are connected by a $p$-closed path in the dual lattice is denoted by $\mathcal{K}_1^* \overset{p^*}{\longleftrightarrow} \mathcal{K}_2^*$.

For positive integers $m < n$, $k$ and $p \in [0, 1]$, let $A_{n,p}$ be the event that there is a $p$-open circuit around the origin of diameter at least $n$ and let $B_{n,p}$ be the event that there is a $p$-closed circuit around the origin in the dual lattice of diameter at least $n$. Let $A_{m,n,p}$ be the event that there is a $p$-open circuit around the origin in the annulus $\text{Ann}(m, n)$ and let $B_{m,n,p}$ be the event that there is a $p$-closed circuit around the origin in the annulus $\text{Ann}(m, n)^*$. Let $A_{m,n,p}^k$ be the event that there are $k$ disjoint $p$-open paths connecting $B(m)$ to $\partial B(n)$.

For $p \in [0, 1]$, we consider a probability space $(\Omega_p, \mathcal{F}_p, \mathbb{P}_p)$, where $\Omega_p = \{0, 1\}^{\mathbb{E}^2}$, $\mathcal{F}_p$ is the $\sigma$-field generated by the finite-dimensional cylinders of $\Omega_p$ and $\mathbb{P}_p$ is a product measure on $(\Omega_p, \mathcal{F}_p)$, $\mathbb{P}_p = \prod_{e \in \mathbb{E}^2} \mu_e$, where $\mu_e$ is given by $\mu_e(\omega_e = 1) = 1 - \mu_e(\omega_e = 0) = p$ for vectors $(\omega_e)_{e \in \mathbb{E}^2} \in \Omega_p$. We say that an edge $e$ is *open* or *occupied* if $\omega_e = 1$, and $e$ is *closed* or *vacant* if $\omega_e = 0$. We say that an edge $e^*$ is *open* or *occupied* if $e$ is open, and it is *closed* or *vacant* if $e$ is closed. The event that two sets of sites $\mathcal{K}_1, \mathcal{K}_2 \subset \mathbb{Z}^2$ are connected by an open path is denoted by $\mathcal{K}_1 \leftrightarrow \mathcal{K}_2$ and the event that two sets of sites $\mathcal{K}_1^*, \mathcal{K}_2^* \subset \mathbb{Z}^2$ are connected by a closed path in the dual lattice is denoted by $\mathcal{K}_1^* \overset{*}{\leftrightarrow} \mathcal{K}_2^*$.

For positive integers $m < n$ and $k$, let $A_n$ be the event that there is an occupied circuit around the origin of diameter at least $n$ and let $B_n$ be the event that there is a vacant circuit around the origin in the dual lattice of diameter at least $n$. Let $A_{m,n}$ be the event that there is an occupied circuit around the origin in the annulus $\text{Ann}(m, n)$ and let $B_{m,n}$ be the event that there is an vacant circuit around the origin in the annulus $\text{Ann}(m, n)^*$. Let $A_{m,n}^k$ be the event that there are $k$ disjoint occupied paths connecting $B(m)$ to $\partial B(n)$.

For two functions $g$ and $h$ from a set $\mathcal{X}$ to $\mathbb{R}$, we write $g(z) \asymp h(z)$ to indicate that $g(z)/h(z)$ is bounded away from 0 and $\infty$, uniformly in $z \in \mathcal{X}$. Throughout this paper, we write "log" for $\log_2$. We also write $\mathbb{P}_{cr}$ for $\mathbb{P}_{p_c}$. All of the constants ($C_i$) in the proofs are strictly positive and finite. Their exact values may be different from proof to proof.

### 1.3. *Main results.*

#### 1.3.1. *Probability for $k$ points in the first pond.*

THEOREM 1.1. *Let $C(0)$ be the cluster of the origin in Bernoulli bond percolation. For any $k > 0$,*

$$(1) \qquad \mathbb{P}(x_1, \ldots, x_k \in \hat{V}_1) \asymp \mathbb{P}_{cr}(x_1, \ldots, x_k \in C(0)), \qquad x_1, \ldots, x_k \in \mathbb{Z}^2.$$

REMARK 1. The lower bound follows from the observation that the $p_c$-open cluster of the origin is a subset of $\hat{V}_1$.

The reader may ask whether there is a universal constant $c$ such that, for all $k \geq 1$ and $x_1, \ldots, x_k \in \mathbb{Z}^2$,

$$\mathbb{P}(x_1, \ldots, x_k \in \hat{V}_1) \leq c \mathbb{P}_{cr}(x_1, \ldots, x_k \in C(0)).$$

In the next theorem, we show that the answer to the above question is negative.

THEOREM 1.2.

$$\lim_{n \to \infty} \frac{\mathbb{P}(B(n) \subset \hat{V}_1)}{\mathbb{P}_{cr}(B(n) \subset C(0))} = \infty.$$

#### 1.3.2. *Ponds and $p_c$-open clusters.* We now state two theorems which say that invasion ponds can contain several large $p_c$-open clusters. Let $K \geq 2, N \geq 1$, and let $\mathcal{U}(m, K, N)$ be the event that the $m$th pond contains at least $K$ disjoint $p_c$-open clusters of size at least $N$.

THEOREM 1.3. *With probability one, there exist infinitely many values of $m$ for which $\mathcal{U}(m, K, N)$ holds.*

THEOREM 1.4. *There exists $\varepsilon > 0$, independent of $N$ but dependent on $K$, such that*

$$\mathbb{P}(\mathcal{U}(1, K, N) | \hat{R}_1 \geq N) \geq \varepsilon,$$

*where $\hat{R}_1$ is the radius of the first pond.*

#### 1.3.3. *Radii of the ponds.* We define $\hat{R}_j$ to be the radius of the graph $G_{i_j}$, that is, $\hat{R}_j = \max\{|x| : x \in G_{i_j}\}$. We refer the reader to Section 1.1 for the definitions of $i_j$ and $G_{i_j}$. In the next theorem, we give the asymptotics for the radii $\hat{R}_j$.

THEOREM 1.5. *For any $k \geq 1$,*

$$(2) \qquad \mathbb{P}(\hat{R}_k \geq n) \asymp (\log n)^{k-1} \mathbb{P}_{cr}(0 \leftrightarrow \partial B(n)).$$

REMARK 2. Let $\{0 \leftrightarrow_k \partial B(n)\}$ be the event that there is a path connecting the origin to the boundary of $B(n)$ such that at most $k$ of its edges are closed. If this event holds, then we say that the origin is connected to $\partial B(n)$ by an open path with $k$ defects. It is a consequence of the Russo–Seymour–Welsh (RSW) theorem (see [17], Proposition 18) that

$$\mathbb{P}_{cr}(0 \leftrightarrow_k \partial B(n)) \asymp (\log n)^k \mathbb{P}_{cr}(0 \leftrightarrow \partial B(n)).$$

Therefore, Theorem 1.5 implies that, for any $k \geq 1$,

$$\mathbb{P}(\hat{R}_k \geq n) \asymp \mathbb{P}_{cr}(0 \leftrightarrow_{k-1} \partial B(n)).$$

REMARK 3. For $k = 1$, the statement (2) follows from Theorem 1 in [18]. Note that in the case $k = 1$, the lower bound immediately follows from the fact that $C(0) \subset \hat{V}_1$, where $C(0)$ is the $p_c$-open cluster of the origin for Bernoulli bond percolation. However, in the case $k \geq 2$, the lower bound is not trivial.

Let $\bar{R}_k$ be the diameter of the $k$th pond, $\bar{R}_k = \max\{|x - y| : x, y \in \hat{V}_k\}$. Note that $(\bar{R}_k)$ are related to $(\hat{R}_k)$ via the simple inequalities $\hat{R}_1 \leq \bar{R}_1 \leq 2\hat{R}_1$ and $\hat{R}_k - \hat{R}_{k-1} - 1 \leq \bar{R}_k \leq 2\hat{R}_k$ for $k \geq 2$. The next theorem immediately follows from Theorem 1.5 and the fact that $\mathbb{P}_{cr}(0 \leftrightarrow \partial B(n)) \asymp \mathbb{P}_{cr}(0 \leftrightarrow \partial B(2n))$.

THEOREM 1.6. *For every $k \geq 1$,*

$$\mathbb{P}(\bar{R}_k \geq n) \asymp (\log n)^{k-1} \mathbb{P}_{cr}(0 \leftrightarrow \partial B(n)).$$

1.3.4. *Mutual singularity of IPC and IIC.* First, we recall the definition of the incipient infinite cluster from [8]. It is shown in [8] that the limit

$$\nu(E) = \lim_{N \to \infty} \mathbb{P}_{cr}(E | 0 \leftrightarrow \partial B(N))$$

exists for any event $E$ that depends on the state of finitely many edges in $\mathbb{E}^2$. The unique extension of $\nu$ to a probability measure on configurations of open and closed edges exists. Under this measure, the open cluster of the origin is a.s. infinite. It is called the *incipient infinite cluster* (IIC). Recall the definition of the IPC $\mathcal{S}$ from Section 1.1. The next statement is [6], Theorem 3.

THEOREM 1.7. *For any finite $\mathcal{K} \subset \mathbb{E}^2$ and $x \in \mathbb{Z}^2$, let $\mathcal{K}(x) = x + \mathcal{K} \subset \mathbb{E}^2$, $E_{\mathcal{K}} = \{\mathcal{K} \subset \mathcal{S}\}$ and $E'_{\mathcal{K}} = \{\mathcal{K} \subset C(0)\}$. Then,*

$$\lim_{|x| \to \infty} \mathbb{P}(E_{\mathcal{K}(x)} | x \in \mathcal{S}) = \nu(E'_{\mathcal{K}}).$$

The above theorem says that, asymptotically, the distribution of invaded edges near $x$ is given by the IIC measure. In this paper, we show that, globally, the IPC measure and the IIC measure are entirely different.

THEOREM 1.8. *The laws of IPC and IIC are mutually singular.*

1.4. *Structure of the paper.* We define the correlation length and state some of its properties in Section 2. We prove Theorem 1.1 in Section 3 and Theorem 1.2 in Section 4. The proofs of Theorems 1.3 and 1.4 are given in Section 5. In Section 6, we prove Theorem 1.5. Theorem 1.8 is proved in Section 7. After Sections 1 and 2, the remainder of the paper may be read in any order. For the notation in Sections 3–7, we refer the reader to Section 1.2.

**2. Correlation length and preliminary results.** In this section, we define the correlation length that will play a crucial role in our proofs. The correlation length was introduced in [3] and further studied in [9].

2.1. *Correlation length.* For positive integers $m, n$ and $p \in (p_c, 1]$, let

$$\sigma(n, m, p) = \mathbb{P}_p(\text{there is an open horizontal crossing of } [0, n] \times [0, m]).$$

Given $\varepsilon > 0$, we define

(3) $$L(p, \varepsilon) = \min\{n : \sigma(n, n, p) \geq 1 - \varepsilon\}.$$

$L(p, \varepsilon)$ is called the *finite-size scaling correlation length* and it is known that $L(p, \varepsilon)$ scales like the usual correlation length (see [9]). It was also shown in [9] that the scaling of $L(p, \varepsilon)$ is independent of $\varepsilon$, given that it is small enough, that is, there exists $\varepsilon_0 > 0$ such that for all $0 < \varepsilon_1, \varepsilon_2 \leq \varepsilon_0$, we have $L(p, \varepsilon_1) \asymp L(p, \varepsilon_2)$. For simplicity, we will write $L(p) = L(p, \varepsilon_0)$ for the entire paper. We also define

$$p_n = \sup\{p : L(p) > n\}.$$

It is easy to see that $L(p) \to \infty$ as $p \to p_c$ and $L(p) = 1$ for $p$ close to 1. In particular, the probability $p_n$ is well defined. It is clear from the definitions of $L(p)$ and $p_n$, and from the RSW theorem, that for positive integers $k$ and $l$, there exists $\delta_{k,l} > 0$ such that for any positive integer $n$ and for all $p \in [p_c, p_n]$,

$$\mathbb{P}_p(\text{there is an open horizontal crossing of } [0, kn] \times [0, ln)) > \delta_{k,l}$$

and

$\mathbb{P}_p$(there is a closed horizontal dual crossing of $([0, kn] \times [0, ln))^*) > \delta_{k,l}$.

By the FKG inequality and a standard gluing argument [5], Section 11.7, we get that, for positive integers $n$ and $k \geq 2$, and for all $p \in [p_c, p_n]$,

$\mathbb{P}_p(\text{Ann}(n, kn)$ contains an open circuit around the origin$) > (\delta_{2k,k-1})^4$

and

$\mathbb{P}_p(\text{Ann}(n, kn)^*$ contains a closed dual circuit around the origin$) > (\delta_{2k,k-1})^4$.

2.2. *Preliminary results.* For any positive $l$, we define $\log^{(0)} l = l$ and $\log^{(j)} l = \log(\log^{(j-1)} l)$ for all $j \geq 1$, provided the right-hand side is well defined. For $l > 10$, let

$$(4) \qquad \log^* l = \min\{j > 0 : \log^{(j)} l \text{ is well defined and } \log^{(j)} l \leq 10\}.$$

Our choice of the constant 10 is quite arbitrary; we could take any other large enough positive number instead of 10. For $l > 10$, let

$$(5) \qquad p_l(j) = \begin{cases} \inf\left\{p > p_c : L(p) \leq \dfrac{l}{C_* \log^{(j)} l}\right\}, & \text{if } j \in (0, \log^* l), \\ p_c, & \text{if } j \geq \log^* l, \\ 1, & \text{if } j = 0. \end{cases}$$

The value of $C_*$ will be chosen later. Note that there exists a universal constant $L_0(C_*) > 10$ such that $p_l(j)$ are well defined if $l > L_0(C_*)$ and nonincreasing in $l$. The last observation follows from the monotonicity of $L(p)$ and the fact that the functions $l/\log^{(j)} l$ are nondecreasing in $l$ for $j \in (0, \log^* l)$ and $l \geq 3$.

We give the following results without proofs:

1. (Reference [6], (2.10).) There exists a universal constant $D_1$ such that, for every $l > L_0(C_*)$ and $j \in (0, \log^* l)$,

$$(6) \qquad\qquad C_* \log^{(j)} l \leq \frac{l}{L(p_l(j))} \leq D_1 C_* \log^{(j)} l.$$

2. (Reference [9], Theorem 2.) There exists a constant $D_2$ such that, for all $p > p_c$,

$$(7) \qquad\quad \theta(p) \leq \mathbb{P}_p[0 \leftrightarrow \partial B(L(p))] \leq D_2 \mathbb{P}_{cr}[0 \leftrightarrow \partial B(L(p))],$$

   where $\theta(p) = \mathbb{P}_p(0 \leftrightarrow \infty)$ is the percolation function for Bernoulli percolation.

3. (Reference [16], Section 4.) There exists a constant $D_3$ such that, for all $n \geq 1$,

$$(8) \qquad\qquad\qquad \mathbb{P}_{p_n}(B(n) \leftrightarrow \infty) \geq D_3.$$

4. (Reference [9], (3.61).) There exists a constant $D_4$ such that, for all positive integers $r \leq s$,

$$\frac{\mathbb{P}_{cr}(0 \leftrightarrow \partial B(s))}{\mathbb{P}_{cr}(0 \leftrightarrow \partial B(r))} \geq D_4 \sqrt{\frac{r}{s}}. \tag{9}$$

5. Recall that $B_n$ is the event that there is a closed circuit around the origin in the dual lattice with diameter at least $n$. There exist positive constants $D_5$ and $D_6$ such that, for all $p > p_c$,

$$\mathbb{P}_p(B_n) \leq D_5 \exp\left\{-D_6 \frac{n}{L(p)}\right\}. \tag{10}$$

This follows from, for example, [6], (2.6) and (2.8) (see also [17], Lemma 39 and Remark 40).

6. (Reference [17], Proposition 34.) Fix $e = \langle (0,0), (1,0) \rangle$ and let $A_n^{2,2}$ be the event that $e_x$ and $e_y$ are connected to $\partial B(n)$ by open paths, and $e_x^*$ and $e_y^*$ are connected to $\partial B(n)^*$ by closed paths. Note that these four paths are disjoint and alternate. Then,

$$(p_n - p_c)n^2 \mathbb{P}_{cr}(A_n^{2,2}) \asymp 1, \qquad n \geq 1. \tag{11}$$

**3. Proof of Theorem 1.1.** Before we prove Theorem 1.1, we give two lemmas that will be used in the proof. To simplify the notation, we write $0 = x_0$. For positive integers $m < n$ and $x \in \mathbb{Z}^2$, we define the event

$$A_{m,n}(x) = \{\text{there is an open circuit in the annulus } \mathrm{Ann}(x; m, n)\}. \tag{12}$$

LEMMA 3.1. *Given a set of vertices $\{x_1, \ldots, x_k\} \in \mathbb{Z}^2$, let $m_i = \min\{|x_i - x_j| : 0 \leq j \leq k, j \neq i\}$, where $x_0 = 0$ and let $m = \min\{m_i : 0 \leq i \leq k\}$. Furthermore, assume $m = m_k$. There then exists a constant $C_1$, independent of $k$, such that for all $p > p_c$, the probability*

$$\mathbb{P}_p(x_1 \leftrightarrow \infty, \ldots, x_k \leftrightarrow \infty, 0 \leftrightarrow \infty)$$

*is bounded from above by*

$$C_1 \mathbb{P}_p(x_k \leftrightarrow \partial B(x_k, m)) \mathbb{P}_p(x_1 \leftrightarrow \infty, \ldots, x_{k-1} \leftrightarrow \infty, 0 \leftrightarrow \infty).$$

PROOF. The statement is trivial if $m \leq 4$, so we assume that $m > 4$. By the RSW theorem, there is a constant $C_2$ independent of $k$ and $m$ such that, for all $p > p_c$, $\mathbb{P}_p(A_{[m/4],[m/2]}(x_k)) \geq 1/C_2$ and hence $1 \leq C_2 \mathbb{P}_p(A_{[m/4],[m/2]}(x_k))$. The FKG inequality gives

$$\begin{aligned}
&\mathbb{P}_p(x_1 \leftrightarrow \infty, \ldots, x_k \leftrightarrow \infty, 0 \leftrightarrow \infty) \\
&\leq C_2 \mathbb{P}_p(A_{[m/4],[m/2]}(x_k)) \mathbb{P}_p(x_1 \leftrightarrow \infty, \ldots, x_k \leftrightarrow \infty, 0 \leftrightarrow \infty) \\
&\leq C_2 \mathbb{P}_p(A_{[m/4],[m/2]}(x_k), x_1 \leftrightarrow \infty, \ldots, x_k \leftrightarrow \infty, 0 \leftrightarrow \infty).
\end{aligned} \tag{13}$$

The event on the right-hand side of (13) implies the following two events:

1. $\{x_k \leftrightarrow \partial B(x_k, [m/4])\}$;
2. $\{x_1 \leftrightarrow \infty, \ldots, x_{k-1} \leftrightarrow \infty, 0 \leftrightarrow \infty \text{ outside } B(x_k, [m/4])\}$.

These two events are independent and therefore the right-hand side of (13) is bounded from above by

$$C_2 \mathbb{P}_p(x_k \leftrightarrow \partial B(x_k, [m/4]))$$
$$\times \mathbb{P}_p(x_1 \leftrightarrow \infty, \ldots, x_{k-1} \leftrightarrow \infty, 0 \leftrightarrow \infty \text{ outside } B(x_k, [m/4]))$$
$$\leq C_2 \mathbb{P}_p(x_k \leftrightarrow \partial B(x_k, [m/4])) \mathbb{P}_p(x_1 \leftrightarrow \infty, \ldots, x_{k-1} \leftrightarrow \infty, 0 \leftrightarrow \infty),$$

where the last inequality follows from monotonicity. Finally, it follows from the FKG inequality, RSW theorem and a standard gluing argument [5], Section 11.7, that $\mathbb{P}_p(x_k \leftrightarrow \partial B(x_k, [m/4])) \asymp \mathbb{P}_p(x_k \leftrightarrow \partial B(x_k, m))$ uniformly in $p > p_c$. □

We recall the definition of the probabilities $(p_n(j))$ in (5). We also recall that these probabilities are well defined if $n > L_0(C_*)$, where $C_*$ is the constant from (5). Later, we choose $C_*$ to be sufficiently large.

LEMMA 3.2.   *Given a set of vertices* $\{x_1, \ldots, x_k\} \in \mathbb{Z}^2$, *let* $n = \max\{|x_i - x_j| : i, j = 0, \ldots, k\}$, *where* $x_0 = 0$. *Furthermore, assume that* $n \geq L_0(C_*)$. *There is then a universal constant* $C_3$ *such that, for all* $j \in (0, \log^* n)$,

$$(14) \qquad \mathbb{P}_{p_n(j)}(x_1 \leftrightarrow \infty, \ldots, x_k \leftrightarrow \infty, 0 \leftrightarrow \infty)$$
$$\leq (C_3 \log^{(j)} n)^{(k+1)/2} \mathbb{P}_{cr}(x_1, \ldots, x_k \in C(0)).$$

PROOF.   We will use induction in $k$. First, we consider the case $k = 1$. To simplify our notation, we write $x_1 = x$. Note that, now, $|x| = n = m$, where $m$ is defined as in Lemma 3.1. From Lemma 3.1, it follows that

$$(15) \qquad \mathbb{P}_{p_n(j)}(x \leftrightarrow \infty, 0 \leftrightarrow \infty) \leq C_1 \theta(p_n(j)) \mathbb{P}_{p_n(j)}(0 \leftrightarrow \partial B(n)).$$

Since $L(p_n(j)) \leq n$, we obtain

$$\mathbb{P}_{p_n(j)}(0 \leftrightarrow \partial B(n)) \leq \mathbb{P}_{p_n(j)}(0 \leftrightarrow \partial B(L(p_n(j)))).$$

Combined with (6), (7) and (9), the above inequality gives

$$C_1 \theta(p_n(j)) \mathbb{P}_{p_n(j)}(0 \leftrightarrow \partial B(n)) \leq C_4 \mathbb{P}_{cr}(0 \leftrightarrow \partial B(L(p_n(j))))^2$$
$$\leq C_5 \frac{n}{L(p_n(j))} \mathbb{P}_{cr}(0 \leftrightarrow \partial B(n))^2 \leq C_6 \log^{(j)} n \mathbb{P}_{cr}(0 \leftrightarrow \partial B(n))^2.$$

The RSW theorem and the gluing argument show (see, e.g., [7], (4)) that

$$(16) \qquad \mathbb{P}_{cr}(0 \leftrightarrow \partial B(n))^2 \leq C_7 \mathbb{P}_{cr}(x \in C(0))$$

for some constant $C_7$. In particular, (14) follows for $k = 1$.

The general case is more involved. We assume that Lemma 3.2 is proved for any set of vertices $\{y_1, \ldots, y_{k-1}\} \in \mathbb{Z}^2$. Then, for a set of vertices $\{x_1, \ldots, x_k\} \in \mathbb{Z}^2$, we define $m$ as in Lemma 3.1 and assume that $m = m_k = \min\{|x_i - x_k| : i < k\}$. We also define $n_1 = \max\{|x_i - x_j| : i, j = 0, \ldots, k-1\}$, with $x_0 = 0$. Then, by the induction hypothesis,

$$(17) \qquad \begin{aligned} &\mathbb{P}_{p_{n_1}(j)}(x_1 \leftrightarrow \infty, \ldots, x_{k-1} \leftrightarrow \infty, 0 \leftrightarrow \infty) \\ &\qquad \leq (C_3 \log^{(j)} n_1)^{k/2} \mathbb{P}_{cr}(x_1, \ldots, x_{k-1} \in C(0)). \end{aligned}$$

Since $n_1 \leq n$ and $m \leq n$, we get $p_n(j) \leq p_m(j)$ and $p_n(j) \leq p_{n_1}(j)$ (see Section 2). Therefore,

$$\begin{aligned} &\mathbb{P}_{p_n(j)}(x_1 \leftrightarrow \infty, \ldots, x_k \leftrightarrow \infty, 0 \leftrightarrow \infty) \\ &\quad \leq C_1 \mathbb{P}_{p_m(j)}(x_k \leftrightarrow \partial B(x_k, m)) \mathbb{P}_{p_{n_1}(j)}(x_1 \leftrightarrow \infty, \ldots, x_{k-1} \leftrightarrow \infty, 0 \leftrightarrow \infty) \\ &\quad \leq C_1 \mathbb{P}_{p_m(j)}(x_k \leftrightarrow \partial B(x_k, m))(C_3 \log^{(j)} n_1)^{k/2} \mathbb{P}_{cr}(x_1, \ldots, x_{k-1} \in C(0)) \\ &\quad \leq (C_8 \log^{(j)} m)^{1/2} \mathbb{P}_{cr}(x_k \leftrightarrow \partial B(x_k, m)) \\ &\qquad \times (C_3 \log^{(j)} n_1)^{k/2} \mathbb{P}_{cr}(x_1, \ldots, x_{k-1} \in C(0)) \\ &\quad \leq C_8^{1/2} C_3^{k/2} (\log^{(j)} n)^{(k+1)/2} \mathbb{P}_{cr}(x_k \leftrightarrow \partial B(x_k, m)) \mathbb{P}_{cr}(x_1, \ldots, x_{k-1} \in C(0)), \end{aligned}$$

where the first inequality follows from Lemma 3.1 and monotonicity, the second inequality follows from (17) and the third inequality follows from (6) and (9). Note that $C_8$ is independent of $k$. It now suffices to show that there is a universal constant $C_9$ such that

$$(18) \qquad \begin{aligned} &\mathbb{P}_{cr}(x_k \leftrightarrow \partial B(x_k, m)) \mathbb{P}_{cr}(x_1, \ldots, x_{k-1} \in C(0)) \\ &\qquad \leq C_9 \mathbb{P}_{cr}(x_1, \ldots, x_k \in C(0)). \end{aligned}$$

Assume that (18) is proved. We can then take $C_3 = \max\{C_6 C_7, C_8 C_9^2\}$. The argument above shows that we can proceed to the next $k$ using this value of $C_3$. We now show (18). We take $x_i$ such that $m = |x_k - x_i|$. Note that this vertex may be the origin. We know that at least one such vertex exists. Recall the definition of events $A_{m,n}(x)$ from (12). By the RSW theorem, there is a constant $C_{10}$ such that $1 \leq C_{10} \mathbb{P}_{cr}(A_{[m/2],m}(x_i); A_{[m/2],m}(x_k))$. Using the FKG inequality, we get

$$\begin{aligned} &\mathbb{P}_{cr}(x_k \leftrightarrow \partial B(x_k, m)) \mathbb{P}_{cr}(x_1, \ldots, x_{k-1} \in C(0)) \\ &\quad \leq C_{10} \mathbb{P}_{cr}(A_{[m/2],m}(x_i); A_{[m/2],m}(x_k)) \mathbb{P}_{cr}(x_k \leftrightarrow \partial B(x_k, m)) \\ &\qquad \times \mathbb{P}_{cr}(x_1, \ldots, x_{k-1} \in C(0)) \\ &\quad \leq C_{10} \mathbb{P}_{cr}(A_{[m/2],m}(x_i); A_{[m/2],m}(x_k); \\ &\qquad\qquad x_k \leftrightarrow \partial B(x_k, m); x_1, \ldots, x_{k-1} \in C(0)). \end{aligned}$$

We show that the event

$$\{A_{[m/2],m}(x_i); A_{[m/2],m}(x_k); x_k \leftrightarrow \partial B(x_k,m); x_1,\dots,x_{k-1} \in C(0)\}$$

implies the event $\{x_i \leftrightarrow x_k; x_1,\dots,x_{k-1} \in C(0)\}$. Indeed, it follows from simple observations:

1. Since the events $\{x_k \leftrightarrow \partial B(x_k,m)\}$ and $A_{[m/2],m}(x_k)$ hold, $x_k$ is connected to the circuit lying in the annulus $\mathrm{Ann}(x_k; [m/2], m)$.
2. Since the distance between $x_i$ and $x_k$ is $m$, the boxes $B(x_i, [m/2]+1)$ and $B(x_k, [m/2]+1)$ intersect. This implies that the circuits in the annuli $\mathrm{Ann}(x_k; [m/2], m)$ and $\mathrm{Ann}(x_i; [m/2], m)$ intersect.
3. Recall that $m$ is the minimal distance in the graph with vertex set $\{0, x_1, \dots, x_k\}$. Since $k \geq 2$ and $\{x_1, \dots, x_{k-1} \in C(0)\}$, there is a vertex $x_j \neq x_k$ (it may be the origin) such that $x_j \notin B(x_i, m-1)$ and $x_j$ is connected to $x_i$. The last observation implies that $x_i$ is connected to the circuit lying in $\mathrm{Ann}(x_i; [m/2], m)$ and hence also to $x_k$.

This proves (18). $\quad\square$

PROOF OF THEOREM 1.1. For $\{x_1, \dots, x_k\} \in \mathbb{Z}^2$, we define, as in Lemma 3.2, $n = \max\{|x_i - x_j| : i, j = 0, \dots, k\}$. If $n < L_0(C_*)$, then $\mathbb{P}_{cr}(x_1, \dots, x_k \in C(0)) > \mathrm{const}(C_*)$. Theorem 1.1 immediately follows since $\mathbb{P}(x_1, \dots, x_k \in \hat{V}_1) \leq 1$. We can therefore assume that $n \geq L_0(C_*)$. In particular, the probabilities $p_n(j)$ are well defined. The rest of the proof is similar to the proof of Theorem 1 in [18]. Recall that $\hat{\tau}_1$ is the value of the outlet of the first pond. We decompose the event $\{x_1, \dots x_k \in \hat{V}_1\}$ according to the value of $\hat{\tau}_1$. We write

$$(19) \quad \mathbb{P}(x_1, \dots, x_k \in \hat{V}_1) = \sum_{j=1}^{\log^* n} \mathbb{P}(x_1, \dots, x_k \in \hat{V}_1, \hat{\tau}_1 \in [p_n(j), p_n(j-1))).$$

Note that, for any $p > p_c$,

(a) if $\hat{\tau}_1 < p$, then any invaded site is in the infinite $p$-open cluster;
(b) if a given set of vertices $\{x_1, \dots, x_k\}$ is in the first pond, $n$ is defined as in Lemma 3.2 and $\hat{\tau}_1 > p$, then there is a $p$-closed circuit around the origin with diameter at least $n$.

We recall the definition of the event

$$B_{n,p} = \{\exists p\text{-closed circuit around } 0 \text{ in the dual with diameter at least } n\}.$$

We conclude that the probability $\mathbb{P}(x_1, \dots, x_k \in \hat{V}_1, \hat{\tau}_1 \in [p_n(j), p_n(j-1)))$ is bounded from above by

$$(20) \quad \mathbb{P}(x_1 \overset{p_n(j-1)}{\longleftrightarrow} \infty, \dots, x_k \overset{p_n(j-1)}{\longleftrightarrow} \infty, 0 \overset{p_n(j-1)}{\longleftrightarrow} \infty; B_{n,p_n(j)}).$$

The FKG inequality implies that the probability (20) is not bigger than

$$(21) \quad \mathbb{P}_{p_n(j-1)}(x_1 \leftrightarrow \infty, \ldots, x_k \leftrightarrow \infty, 0 \leftrightarrow \infty)\mathbb{P}(B_{n,p_n(j)})$$

$$\leq C_{11}(\log^{(j-1)} n)^{-C_{12}}\mathbb{P}_{p_n(j-1)}(x_1 \leftrightarrow \infty, \ldots, x_k \leftrightarrow \infty, 0 \leftrightarrow \infty),$$

where we use (6) and (10) to bound the probability of $B_{n,p_n(j)}$ by $C_{11}(\log^{(j-1)} n)^{-C_{12}}$. The constant $C_{12}$ can be made arbitrarily large provided that $C_*$ is made large enough. We consider bounds for (21) separately for $j = 1$ and for $j > 1$. If $j > 1$, we use Lemma 3.2 to bound (21) by

$$C_{13}(\log^{(j-1)} n)^{(k+1)/2 - C_{12}}\mathbb{P}_{cr}(x_1, \ldots, x_k \in C(0)).$$

If $j = 1$, we bound (21) by

$$C_{11}n^{-C_{12}} \leq C_{14}n^{-1/2}\mathbb{P}_{cr}(0 \leftrightarrow \partial B(n))^{2k} \leq C_{15}n^{-1/2}\mathbb{P}_{cr}(x_1, \ldots, x_k \in C(0)).$$

The first inequality holds for $C_{12} \geq k + 1/2$ since $\mathbb{P}_{cr}(0 \leftrightarrow \partial B(n)) > \frac{1}{2}n^{-1/2}$ (see [5], (11.90)). The last inequality follows from (16), applied $k$ times, and the FKG inequality. Therefore, for all $j$, if $C_{12} \geq k + 1/2$, then (21) is bounded by

$$C_{16}(\log^{(j-1)} n)^{-1/2}\mathbb{P}_{cr}(x_1, \ldots, x_k \in C(0)).$$

We plug this bound into (19):

$$\mathbb{P}(x_1, \ldots, x_k \in \hat{V}_1) \leq C_{16}\mathbb{P}_{cr}(x_1, \ldots, x_k \in C(0)) \sum_{j=1}^{\log^* n} (\log^{(j-1)} n)^{-1/2}$$

$$\leq C_{17}\mathbb{P}_{cr}(x_1, \ldots, x_k \in C(0)).$$

The last inequality follows from the fact that

$$\sup_{n>10} \sum_{j=1}^{\log^* n} (\log^{(j-1)} n)^{-1/2} < \infty$$

(see, e.g., [6], (2.26)).  □

## 4. Proof of Theorem 1.2.

In this section, we prove that

$$\lim_{n \to \infty} \frac{\mathbb{P}(B(n) \subset \hat{V}_1)}{\mathbb{P}_{cr}(B(n) \subset C(0))} = \infty.$$

By RSW arguments [5], Section 11.7, the denominator is at most equal to

$$C_1\mathbb{P}_{cr}(B(n) \subset C(0) \text{ in } B(2n))$$

for some $C_1 > 0$. Recall that $p_n = \sup\{p : L(p) > n\}$. We can bound the numerator from below: it is at least equal to

$$\mathbb{P}_{p_n}(B(n) \subset C(0) \text{ in } B(2n) \cap \exists \text{ closed circuit around } B(2n))$$

$$= \mathbb{P}_{p_n}(B(n) \subset C(0) \text{ in } B(2n))\mathbb{P}_{p_n}(\exists \text{ closed circuit around } B(2n)).$$

By the definition of $L(p)$, there exists $C_2 > 0$ such that this probability is at least

$$C_2\mathbb{P}_{p_n}(B(n) \subset C(0) \text{ in } B(2n)).$$

Therefore, to prove Theorem 1.2, it suffices to show that

$$(22) \qquad \lim_{n\to\infty} \frac{\mathbb{P}_{p_n}(B(n) \subset C(0) \text{ in } B(2n))}{\mathbb{P}_{cr}(B(n) \subset C(0) \text{ in } B(2n))} = \infty.$$

For this, we use Russo's formula [5] (the definition of pivotal edges is also given in [5]). Let $\Gamma_n$ be the event which appears both in the numerator and in the denominator of (22). Let $p \in [\epsilon, 1 - \epsilon]$ for some $\epsilon < \frac{1}{2}$ and, for any vertex $v$, let $E_v$ be the set of edges incident to $v$. We see that

$$\frac{d}{dp}\mathbb{P}_p(\Gamma_n) = \sum_e \mathbb{P}_p(e \text{ is pivotal for } \Gamma_n)$$

$$\geq \frac{1}{2p} \sum_{v \in B(n)} \sum_{e \in E_v} \mathbb{P}_p(e \text{ is pivotal for } \Gamma_n; \Gamma_n)$$

$$\geq \frac{1}{2p} \sum_{v \in B(n)} \mathbb{P}_p(\exists e \in E_v \text{ pivotal for } \Gamma_n; \Gamma_n)$$

$$\geq \frac{1}{2p} \sum_{v \in B(n)} \min(p, 1-p)^4 \mathbb{P}_p(\Gamma_n)$$

$$\geq C_3 n^2 \mathbb{P}_p(\Gamma_n).$$

In particular,

$$\mathbb{P}_{p_n}(\Gamma_n) \geq \mathbb{P}_{cr}(\Gamma_n)e^{C_4 n^2(p_n - p_c)}$$

for some $C_4 > 0$. It follows from (11) and the fact that $\theta(p_c) = 0$ that $n^2(p_n - p_c) \to \infty$. This completes the proof.

**5. Proofs of Theorems 1.3 and 1.4.** First, we prove two lemmas (see Section 1.2 for the definitions).

LEMMA 5.1. *For each $k \geq 2$, there exists $c_k$ such that, for all $n$,*

$$\mathbb{P}(A^1_{n,kn,p_c}) \leq c_k,$$

*where $c_k \to 0$ as $k \to \infty$.*

PROOF. Recall that $B_{n,2n} = \{$there is a closed circuit in $\text{Ann}(n,2n)^*\}$. Pick $c > 0$ such that, for all $N \geq 1$,

$$\mathbb{P}_{cr}(B_{N,2N}) \geq c.$$

We split the annulus $\text{Ann}(n,kn)$ into $[\log k]$ disjoint annuli $\text{Ann}(2^i n, 2^{i+1} n)$:

$$\mathbb{P}(A^1_{n,kn,p_c}) \leq (1-c)^{\log k - 1}.$$

This completes the proof. $\square$

LEMMA 5.2. *There exists $C_1 > 0$ such that for all $N$ and $k$,*

$$\mathbb{P}(A_{N,2N,p_c} \cap A_{kN,2kN,p_c} \cap A^1_{N,2kN,p_N}) \geq C_1.$$

PROOF. By RSW arguments, there exists $C_2 > 0$ such that for all $N$ and $k$,

$$\mathbb{P}(A_{N,2N,p_c} \cap A_{kN,2kN,p_c}) \geq C_2.$$

If follows from (8) that there exists $C_3 > 0$ such that for all $N$ and $k$,

$$\mathbb{P}(A^1_{N,2kN,p_N}) \geq \mathbb{P}_{p_N}(B(N) \leftrightarrow \infty) \geq C_3.$$

The FKG inequality gives the result. $\square$

We now prove the theorems.

PROOF OF THEOREM 1.3. We prove the theorem for $K = 2$. For other values of $K$, the proof is similar. Let $D(k,N) = A_{N,2N,p_c} \cap A_{kN,2kN,p_c} \cap A^1_{N,2kN,p_N}$ and pick $C_1$ from Lemma 5.2. Fix $k$ such that the constant $c_{k/2}$ from Lemma 5.1 satisfies $c_{k/2} \leq \frac{C_1}{2}$. It follows that

$$\mathbb{P}(D(k,N) \cap \{A^1_{2N,kN,p_c}\}^c) \geq \frac{C_1}{2}.$$

For any $k \geq 2$, there exists $C_4 = C_4(k)$ such that for all $N$,

$$\mathbb{P}(B_{2kN,4kN,p_N}) \geq C_4.$$

Therefore, by independence,

$$\mathbb{P}(D(k,N) \cap \{A^1_{2N,kN,p_c}\}^c \cap B_{2kN,4kN,p_N}) \geq \frac{C_1 C_4}{2} > 0.$$

This statement, along with the Borel–Cantelli lemma, gives the theorem. $\square$

PROOF OF THEOREM 1.4. Let $A^1_{n,p} = \{0 \leftrightarrow \partial B(n) \text{ by a } p\text{-open path}\}$. We first note that [18] gives a constant $C_5 > 0$ such that for all $N$,

$$\mathbb{P}(\hat{R}_1 \geq N) \leq C_5 \mathbb{P}(A^1_{2N,p_c}).$$

It is obvious that $\mathbb{P}(\hat{R}_1 \geq N \cap \mathcal{U}(1,2,N)) \geq \mathbb{P}(A^1_{2N,p_c} \cap \mathcal{U}(1,2,N))$. Therefore, it suffices to show that there is an $\varepsilon > 0$ such that for all $N$,

$$\mathbb{P}(\mathcal{U}(1,2,N)|A^1_{2N,p_c}) \geq \varepsilon.$$

The rest of the proof is almost the same as the proof of Theorem 1.3. Let $D(k,N)$ be as in the proof of Theorem 1.3. Pick $C_1$ from Lemma 5.2. By the FKG inequality, we see that

$$\mathbb{P}(D(k,N) \cap A^1_{2N,p_c}) \geq C_1 \mathbb{P}(A^1_{2N,p_c}).$$

By independence and Lemma 5.1, we may fix $k$ such that for all $N$,

$$\mathbb{P}(A^1_{2N,p_c} \cap A^1_{2N,kN,p_c}) \leq c_{k/2} \mathbb{P}(A^1_{2N,p_c}) \leq \frac{C_1}{2} \mathbb{P}(A^1_{2N,p_c}).$$

For any $k \geq 2$, there exists $C_4 = C_4(k)$ such that for all $N$,

$$\mathbb{P}(B_{2kN,4kN,p_N}) \geq C_4.$$

Independence now gives us

$$\mathbb{P}(A^1_{2N,p_c} \cap D(k,N) \cap \{A^1_{2N,kN,p_c}\}^c \cap B_{2kN,4kN,p_N}) \geq \frac{C_1 C_4}{2} \mathbb{P}(A^1_{2N,p_c}).$$

This concludes the proof.  □

## 6. Proof of Theorem 1.5.

6.1. *Upper bound.* We give the proof for $k = 2$. The case $k = 1$ is considered in [18] and the proof for $k \geq 3$ is similar to the proof for $k = 2$.

We fix $n$ and divide the box $B(n)$ into $[\log n] + 1$ annuli. We write

$$(23) \quad \mathbb{P}(\hat{R}_2 \geq n) = \mathbb{P}(\hat{R}_1 \geq n) + \sum_{k=1}^{[\log n]+1} \mathbb{P}\left(\hat{R}_2 \geq n, \hat{R}_1 \in \left[\frac{n}{2^k}, \frac{n}{2^{k-1}}\right)\right).$$

Since [18], Theorem 1,

$$\mathbb{P}(\hat{R}_1 \geq n) \leq C_1 \mathbb{P}_{cr}(0 \leftrightarrow \partial B(n)) \leq C_1 \log n \mathbb{P}_{cr}(0 \leftrightarrow \partial B(n)),$$

it remains to bound the typical term of the sum on the right-hand side of (23). It is sufficient to show that there exists a constant $C_2$ such that, for any $m \in [0, n/2]$,

$$(24) \quad \mathbb{P}(\hat{R}_2 \geq n; \hat{R}_1 \in [m, 2m]) \leq C_2 \mathbb{P}_{cr}(0 \leftrightarrow \partial B(n)).$$

We only consider the case $m \geq L_0(C_*)$. The proof for $m < L_0(C_*)$ is similar to the proof for $m \geq L_0(C_*)$, but much simpler. We omit the details. We now assume that $m \geq L_0(C_*)$. In particular, the probabilities $(p_m(i))$ and $(p_n(j))$ are well defined.

We decompose the event on the left-hand side according to the $\tau$ value of the first and the second outlet. The probability $\mathbb{P}(\hat{R}_2 \geq n; \hat{R}_1 \in [m, 2m])$ is bounded from above by

$$
\begin{aligned}
(25) \quad \sum_{i=1}^{\log^* m} \sum_{j=1}^{\log^* n} \mathbb{P}(\hat{R}_2 \geq n; \hat{R}_1 \in [m, 2m]; \\
\hat{\tau}_1 \in [p_m(i), p_m(i-1)]; \hat{\tau}_2 \in [p_n(j), p_n(j-1)]).
\end{aligned}
$$

Note that if the event $\{\hat{R}_1 \geq m; \hat{\tau}_1 \in [p_m(i), p_m(i-1)]\}$ occurs, then:

- there is a $p_m(i-1)$-open path from the origin to infinity;
- the origin is surrounded by a $p_m(i)$-closed circuit of diameter at least $m$ in the dual lattice.

We also note that if the event $\{\hat{R}_1 \leq 2m; \hat{R}_2 \geq n; \hat{\tau}_2 \in [p_n(j), p_n(j-1)]\}$ occurs, then:

- there is a $p_n(j-1)$-open path from the box $B(2m)$ to infinity;
- the origin is surrounded by a $p_n(j)$-closed circuit of diameter at least $n$ in the dual lattice.

From the two observations above, the sum (25) is less than

$$
(26) \quad \sum_{i=1}^{\log^* m} \sum_{j=1}^{\log^* n} \mathbb{P}(0 \overset{p_m(i-1)}{\longleftrightarrow} \partial B(m); B(2m) \overset{p_n(j-1)}{\longleftrightarrow} \partial B(n); B_{m,p_m(i)}; B_{n,p_n(j)}).
$$

The FKG inequality and the independence of the first two events together imply that (26) is not larger than

$$
\begin{aligned}
(27) \quad \sum_{i=1}^{\log^* m} \sum_{j=1}^{\log^* n} \mathbb{P}_{p_m(i-1)}(0 \leftrightarrow \partial B(m)) \mathbb{P}_{p_n(j-1)}(B(2m) \leftrightarrow \partial B(n)) \\
\times \mathbb{P}(B_{m,p_m(i)}; B_{n,p_n(j)}).
\end{aligned}
$$

We use (6) and (10) to bound the probability of $B_{m,p_m(i)}$ by $C_3(\log^{(i-1)} m)^{-C_4}$, where $C_4$ can be made arbitrarily large, provided that $C_*$ is made large enough. Substitution gives a bound for the last term of (27):

$$
\begin{aligned}
(28) \quad \mathbb{P}(B_{m,p_m(i)}; B_{n,p_n(j)}) &\leq \min[C_3(\log^{(i-1)} m)^{-C_4}, C_3(\log^{(j-1)} n)^{-C_4}] \\
&= C_3 \max[\log^{(i-1)} m, \log^{(j-1)} n]^{-C_4} \\
&\leq C_3(\log^{(i-1)} m)^{-C_4/2}(\log^{(j-1)} n)^{-C_4/2}.
\end{aligned}
$$

The RSW theorem and the FKG inequality together imply that

$$
(29) \quad \mathbb{P}_p(0 \leftrightarrow \partial B(m)) \mathbb{P}_p(B(2m) \leftrightarrow \partial B(n)) \leq C_5 \mathbb{P}_p(0 \leftrightarrow \partial B(n)),
$$

uniformly in $p \geq p_c$. Furthermore, using (6)–(9), we get

$$(30) \qquad \mathbb{P}_{p_m(i-1)}(0 \leftrightarrow \partial B(m)) \leq C_6 (\log^{(i-1)} m)^{1/2} \mathbb{P}_{cr}(0 \leftrightarrow \partial B(m))$$

and

$$(31) \qquad \mathbb{P}_{p_n(j-1)}(B(2m) \leftrightarrow \partial B(n)) \leq C_7 (\log^{(j-1)} n)^{1/2} \mathbb{P}_{cr}(B(2m) \leftrightarrow \partial B(n)).$$

In the last inequality, we also use (29). We apply the inequalities (28), (29), (30) and (31) to (27). We obtain that the probability $\mathbb{P}(\hat{R}_2 \geq n; \hat{R}_1 \in [m, 2m])$ is not larger than

$$C_8 \mathbb{P}_{cr}(0 \leftrightarrow \partial B(n)) \sum_{i=1}^{\log^* m} \sum_{j=1}^{\log^* n} (\log^{(i-1)} m)^{-(C_4-1)/2} (\log^{(j-1)} n)^{-(C_4-1)/2}.$$

We take $C_*$ large enough so that $C_4$ is greater than 1. As in (2.26) of [6], it is easy to see that there exists a universal constant $C_9 < \infty$ such that for all $n > 10$,

$$\sum_{j=1}^{\log^* n} (\log^{(j-1)} n)^{-(C_4-1)/2} \leq C_9.$$

6.2. *Lower bound.* We first give the main idea of the proof. Recall from Remark 2 that it is equivalent to prove that $\mathbb{P}(\hat{R}_k \geq n) \geq c_k \mathbb{P}_{cr}(0 \leftrightarrow_{k-1} \partial B(n))$ for some positive constants $c_k$ that do not depend on $n$. In the case $k = 1$, the event $\{0 \xleftrightarrow{p_c} \partial B(n)\}$ obviously implies the event $\{\hat{R}_1 \geq n\}$. However, for $k \geq 2$, the event $\{0 \xleftrightarrow{p_c}_{k-1} \partial B(n)\}$ does not, in general, imply the event $\{\hat{R}_k \geq n\}$. The weights of some defected edges from the definition of the event $\{0 \xleftrightarrow{p_c}_{k-1} \partial B(n)\}$ can be large enough so that these edges are never invaded. We resolve this problem by constructing a subevent of the event $\{0 \xleftrightarrow{p_c}_{k-1} \partial B(n)\}$ which implies the event $\{\hat{R}_k \geq n\}$ and, moreover, the probability of this new event is comparable with the probability $\mathbb{P}(0 \xleftrightarrow{p_c}_{k-1} \partial B(n))$. To construct such an event, we first extend results from [9] in Lemmas 6.2 and 6.3 below. We then construct events that will be used in the proof of the lower bound in Theorem 1.5 and show that they satisfy the desired properties (see, e.g., Corollary 6.2 below).

We begin with some definitions and lemmas.

LEMMA 6.1 (Generalized FKG). *Let $\xi_1, \ldots, \xi_n$ be i.i.d. real-valued random variables. Let $I_1, I_2, I_3$ be disjoint subsets of $\{1, \ldots, n\}$. Let $A_1 \in \sigma(\xi_i : i \in I_1 \cup I_2)$ and $A_2 \in \sigma(\xi_i : i \in I_2)$ be increasing in $(\xi_i)$. Let $B_1 \in \sigma(\xi_i : i \in I_1 \cup I_3)$ and $B_2 \in \sigma(\xi_i : i \in I_3)$ be decreasing in $(\xi_i)$. Then,*

$$(32) \qquad \mathbb{P}(A_2 \cap B_2 | A_1 \cap B_1) \geq \mathbb{P}(A_2)\mathbb{P}(B_2).$$

Proof. Inequality (32) for $\mathbb{P}_p$ (rather than $\mathbb{P}$) is given in [9], Lemma 3, or [17], Lemma 13. The main ingredient of that proof is the Harris–FKG inequality for $\mathbb{P}_p$ (see [5], Theorem 2.4), which is also valid for $\mathbb{P}$ (see, e.g., [12], Theorem 5.13). Apart from that, the proof of (32) is analogous to the proofs of [9], Lemma 3, and [17], Lemma 13, and so we omit it. □

Although we will not apply Lemma 6.1 to the following events, they serve as simple examples. The events $\{0 \xleftrightarrow{p} \partial B(n)\}$, $\{B(m) \xleftrightarrow{p} \partial B(n)\}$ are decreasing in $(\tau_e)$ and the events $\{0^* \xleftrightarrow{p^*} \partial B(n)^*\}$, $\{B(m)^* \xleftrightarrow{p^*} \partial B(n)^*\}$ are increasing in $(\tau_e)$.

Recall that the ends of an edge $e \in \mathbb{E}^2$, left (resp., right) or bottom (resp., top), are denoted by $e_x, e_y \in \mathbb{Z}^2$ and the ends of its dual edge $e^*$, bottom (resp., top) or left (resp., right), are denoted by $e_x^*$ and $e_y^*$. We also write $(1,0)$ for the edge with ends $(0,0),(1,0) \in \mathbb{Z}^2$.

Definition 6.1. For any positive integer $n$, $q_1, q_2 \in [0,1]$, $z \in \mathbb{Z}^2$ and an edge $e \in B(z,n)$, we define $A_e(z;n;q_1,q_2)$ as the event that there exist four disjoint paths $P_1$–$P_4$ such that:

- $P_1$ and $P_2$ are $q_1$-open paths in $B(z,n) \setminus \{e\}$, the path $P_1$ connects $e_x$ to $\partial B(z,n)$ and the path $P_2$ connects $e_y$ to $\partial B(z,n)$;
- $P_3$ and $P_4$ are $q_2$-closed paths in $B(z,n)^* \setminus \{e^*\}$, the path $P_3$ connects $e_x^*$ to $\partial B(z,n)^*$ and the path $P_4$ connects $e_y^*$ to $\partial B(z,n)^*$.

We write $A_e(n;q_1,q_2)$ for $A_e(0;n;q_1,q_2)$ and $A(n;q_1,q_2)$ for $A_{(1,0)}(n;q_1,q_2)$.

For any two positive integers $n < N$, $q_1, q_2 \in [0,1]$, $z \in \mathbb{Z}^2$, we define $A(z;n,N;q_1,q_2)$ as the event that there exist four disjoint paths, two $q_1$-open paths in the annulus $\text{Ann}(z;n,N)$ from $B(z,n)$ to $\partial B(z,N)$ and two $q_2$-closed paths in the annulus $\text{Ann}(z;n,N)^*$ from $B(z,n)^*$ to $\partial B(z,N)^*$, such that the $q_1$-open paths are separated by the $q_2$-closed paths. We write $A(n,N;q_1,q_2)$ for $A(0;n,N;q_1,q_2)$. The events $A_e(n;q_1,q_2)$ and $A(n,N;q_1,q_2)$ are illustrated in Figure 1.

We will follow the ideas developed in [9]. For that, we need to define some subevents of $A_e(z;n;q_1,q_2)$ and $A(z;n,N;q_1,q_2)$. For $n \geq 1$, let $U_n = \partial B(n) \cap \{x_2 = n\}$, $D_n = \partial B(n) \cap \{x_2 = -n\}$, $R_n = \partial B(n) \cap \{x_1 = n\}$ and $L_n = \partial B(n) \cap \{x_1 = -n\}$ be the sides of the box $B(n)$. Let $U_n(z) = z + U_n$, $D_n(z) = z + D_n$, $R_n(z) = z + R_n$ and $L_n(z) = z + L_n$ be the sides of the box $B(z,n)$.

Definition 6.2. For any positive integer $n$, $q_1, q_2 \in [0,1]$, $z \in \mathbb{Z}^2$ and an edge $e \in B(z,n)$, we define $\bar{A}_e(z;n;q_1,q_2)$ as the event that there exist four disjoint paths $P_1$–$P_4$ such that:

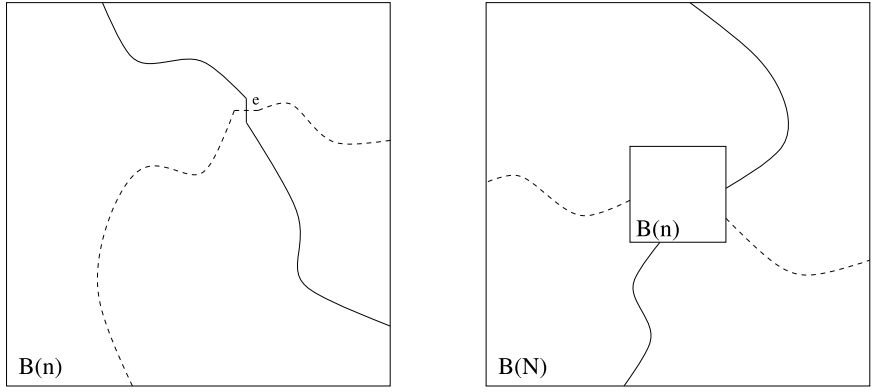

Fig. 1.   *Events $A_e(n; q_1, q_2)$ and $A(n, N; q_1, q_2)$. The solid curves represent $q_1$-open paths, and the dotted curves represent $q_2$-closed paths. The edge $e$ does not have to be $q_1$-open or $q_2$-closed.*

- $P_1$ and $P_2$ are $q_1$-open paths in $B(z, n) \setminus \{e\}$, the path $P_1$ connects $e_x$ or $e_y$ to $U_n(z)$ and the path $P_2$ connects the other end of $e$ to $D_n(z)$;
- $P_3$ and $P_4$ are $q_2$-closed paths in $B(z, n)^* \setminus \{e^*\}$, the path $P_3$ connects $e_x^*$ or $e_y^*$ to $R_n(z)^*$ and the path $P_4$ connects the other end of $e^*$ to $L_n(z)^*$.

We define $\bar{A}_e(z; n; q_1, \cdot)$ as the event that there exist two disjoint $q_1$-open paths $P_1$ and $P_2$ in $B(z, n) \setminus \{e\}$, the path $P_1$ connects $e_x$ or $e_y$ to $U_n(z)$ and the path $P_2$ connects the other end of $e$ to $D_n(z)$.

We write $\bar{A}_e(n; q_1, q_2)$ for $\bar{A}_e(0; n; q_1, q_2)$, $\bar{A}(n; q_1, q_2)$ for $\bar{A}_{(1,0)}(n; q_1, q_2)$ and we use similar notation for the events $A_e(z; n; q_1, \cdot)$.

For any two positive integers $n < N$, $q_1, q_2 \in [0, 1]$ and $z \in \mathbb{Z}^2$, we define $\bar{A}(z; n, N; q_1, q_2)$ as the event that there exist four disjoint paths $P_1$–$P_4$ such that:

- $P_1$ and $P_2$ are $q_1$-open paths in the annulus $\operatorname{Ann}(z; n, N)$, the path $P_1$ connects $U_n(z)$ to $U_N(z)$ and the path $P_2$ connects $D_n(z)$ to $D_N(z)$;
- $P_3$ and $P_4$ are $q_2$-closed paths in the annulus $\operatorname{Ann}(z; n, N)^*$, the path $P_3$ connects $R_n(z)^*$ to $R_N(z)^*$ and the path $P_4$ connects $L_n(z)^*$ to $L_N(z)^*$.

We write $\bar{A}(n, N; q_1, q_2)$ for $\bar{A}(0; n, N; q_1, q_2)$.

We also need to define events similar to the events $\Delta$ in [9], Figure 8. For any two positive integers $n < N$ and $z \in \mathbb{Z}^2$, we define $U_{n,N}(z) = z + [-n, n] \times [n + 1, N]$, $D_{n,N}(z) = z + [-n, n] \times [-N, -n - 1]$, $R_{n,N}(z) = z + [n + 1, N] \times [-n, n]$ and $L_{n,N}(z) = z + [-N, -n - 1] \times [-n, n]$.

DEFINITION 6.3.   For any positive integer $n$, $q_1, q_2 \in [0, 1]$, $z \in \mathbb{Z}^2$ and an edge $e \in B(z, [n/2])$, we define $\bar{A}_e(z; n; q_1, q_2)$ as the event that:

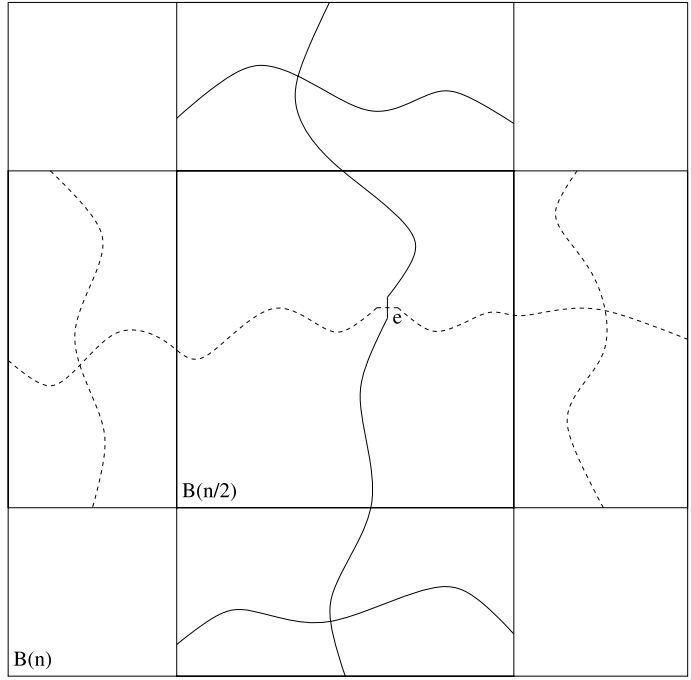

Fig. 2. *Event $\widetilde{A}_e(n; q_1, q_2)$. The solid curves represent $q_1$-open paths and the dotted curves represent $q_2$-closed paths. The edge $e$ does not have to be $q_1$-open or $q_2$-closed.*

- the event $\bar{A}_e(z; n; q_1, q_2)$ occurs;
- the two $q_1$-open paths $P_1$ and $P_2$ from the definition of $\bar{A}_e(z; n; q_1, q_2)$ satisfy $P_1 \cap \mathrm{Ann}(z; [n/2], n) \subset U_{[n/2],n}(z)$ and $P_2 \cap \mathrm{Ann}(z; [n/2], n) \subset D_{[n/2],n}(z)$;
- the two $q_2$-closed paths $P_3$ and $P_4$ from the definition of $\bar{A}_e(z; n; q_1, q_2)$ satisfy $P_3 \cap \mathrm{Ann}(z; [n/2], n)^* \subset R_{[n/2],n}(z)^*$ and $P_4 \cap \mathrm{Ann}(z; [n/2], n)^* \subset L_{[n/2],n}(z)^*$;
- there exist $q_1$-open horizontal crossings of $U_{[n/2],n}(z)$ and $D_{[n/2],n}(z)$ and there exist $q_2$-closed vertical crossings of $L_{[n/2],n}(z)^*$ and $R_{[n/2],n}(z)^*$.

We write $\widetilde{A}_e(n; q_1, q_2)$ for $\widetilde{A}_e(0; n; q_1, q_2)$ and $\widetilde{A}(n; q_1, q_2)$ for $\widetilde{A}_{(1,0)}(n; q_1, q_2)$. The event $\widetilde{A}_e(n; q_1, q_2)$ is illustrated in Figure 2.

For any positive integers $n, N$ such that $4n \leq N$, $q_1, q_2 \in [0, 1]$, $z \in \mathbb{Z}^2$, we define $\widetilde{A}(z; n, N; q_1, q_2)$ as the event that:

- the event $\bar{A}(z; n, N; q_1, q_2)$ occurs;
- the two $q_1$-open paths $P_1$ and $P_2$ from the definition of $\bar{A}(z; n, N; q_1, q_2)$ satisfy $P_1 \cap \mathrm{Ann}(z; n, 2n) \subset U_{n,2n}(z)$, $P_1 \cap \mathrm{Ann}(z; [N/2], N) \subset U_{[N/2],N}(z)$, $P_2 \cap \mathrm{Ann}(z; n, 2n) \subset D_{n,2n}(z)$ and $P_2 \cap \mathrm{Ann}(z; [N/2], N) \subset D_{[N/2],N}(z)$;

- the two $q_2$-closed paths $P_3$ and $P_4$ from the definition of $\bar{A}(z; n, N; q_1, q_2)$ satisfy $P_3 \cap \text{Ann}(z; n, 2n)^* \subset R_{n,2n}(z)^*$, $P_3 \cap \text{Ann}(z; [N/2], N)^* \subset R_{[N/2],N}^*(z)$, $P_4 \cap \text{Ann}(z; n, 2n)^* \subset L_{n,2n}(z)^*$ and $P_4 \cap \text{Ann}(z; [N/2], N)^* \subset L_{[N/2],N}(z)^*$;
- there exist $q_1$-open horizontal crossings of $U_{n,2n}(z)$, $U_{[N/2],N}(z)$, $D_{n,2n}(z)$ and $D_{[N/2],N}(z)$, and there exist $q_2$-closed vertical crossings of $L_{n,2n}(z)^*$, $L_{[N/2],N}(z)^*$, $R_{n,2n}(z)^*$ and $R_{[N/2],N}(z)^*$.

We write $\tilde{A}(n, N; q_1, q_2)$ for $\tilde{A}(0; n, N; q_1, q_2)$.

LEMMA 6.2. *For any positive integers $n, N$ such that $4n \leq N$ and $q_1, q_2 \in [p_c, p_N]$,*

$$(33) \qquad \mathbb{P}(A(n, N; q_1, q_2)) \asymp \mathbb{P}(\bar{A}(n, N; q_1, q_2)) \asymp \mathbb{P}(\tilde{A}(n, N; q_1, q_2))$$

*and*

$$(34) \qquad \mathbb{P}(A(N; q_1, q_2)) \asymp \mathbb{P}(\bar{A}(N; q_1, q_2)) \asymp \mathbb{P}(\tilde{A}(N; q_1, q_2)),$$

*where the constants in (33) and (34) do not depend on $n$, $N$, $q_1$ and $q_2$.*

PROOF. The case $q_1 = q_2$ is considered in [9], Lemma 4 (see also [17], Theorem 11). The proof is based on Lemma 6.1 and the RSW theorem. The same proof is valid for general $q_1$ and $q_2$. □

We need several corollaries of Lemmas 6.1 and 6.2. Their proofs are similar to the proofs for $q_1 = q_2$ (see, e.g., Corollary 3 and Lemma 6 in [9] or Propositions 12 and 17 in [17]). We omit the details.

COROLLARY 6.1. 1. *For any positive integers $a$, $b$ and $n < N$ such that $an < bN$, for any $q_1, q_2 \in [p_c, p_N]$,*

$$(35) \qquad \mathbb{P}(A(n, N; q_1, q_2)) \asymp \mathbb{P}(A(an, bN; q_1, q_2)),$$

*where the constants in (35) only depend on $a$ and $b$.*

2. *For any positive integers $n < m < N$ and $q_1, q_2 \in [p_c, p_N]$,*

$$(36) \qquad \mathbb{P}(A(n, N; q_1, q_2)) \asymp \mathbb{P}(A(n, m; q_1, q_2))\mathbb{P}(A(m, N; q_1, q_2)),$$

*where the constants in (36) do not depend on $n$, $m$, $N$, $q_1$ and $q_2$.*

3. *For any positive integer $N$, $q_1, q_2 \in [p_c, p_N]$ and edge $e \in B([N/2])$,*

$$(37) \qquad \begin{aligned} \mathbb{P}(A_e(N; q_1, q_2)) &\asymp \mathbb{P}(\bar{A}_e(N; q_1, q_2)) \asymp \mathbb{P}(\tilde{A}_e(N; q_1, q_2)) \\ &\asymp \mathbb{P}(A(N; q_1, q_2)), \end{aligned}$$

*where the constants in (37) do not depend on $N$, $q_1$, $q_2$ and $e$.*

The proof of the lower bound in Theorem 1.5 is based on the following lemma.

LEMMA 6.3. *For any positive integer $N$, $q_1, q_2 \in [p_c, p_N]$ and $e \in B([N/2])$,*

$$\mathbb{P}(A_e(N; q_1, q_2)) \asymp \mathbb{P}(A(N; p_c, p_c)), \tag{38}$$

*where the constants in (38) do not depend on $N$, $q_1$, $q_2$ and $e$.*

PROOF. The proof for $q_1 = q_2$ is given in [9], Lemma 8, and [17], Theorem 27. In this case, the probability measure $\mathbb{P}$ can be replaced by the probability measure $\mathbb{P}_{q_1}$ on configurations of open and closed edges. This is not the case when $q_1 \neq q_2$, which makes the proof of (38) more involved. Note that, by (34) and (37), it is sufficient to show that, for $q_1, q_2 \in [p_c, p_N]$,

$$\mathbb{P}(\bar{A}(N; q_1, q_2)) \asymp \mathbb{P}(\bar{A}(N; p_c, p_c)).$$

It is immediate from monotonicity in $q_1$ and $q_2$ that

$$\mathbb{P}(\bar{A}(N; p_c, q_2)) \leq \mathbb{P}(\bar{A}(N; q_1, q_2)) \leq \mathbb{P}(\bar{A}(N; q_1, p_c)).$$

Therefore, it remains to show that there exist constants $D_1$ and $D_2$ such that for all $q_1, q_2 \in [p_c, p_N]$,

$$\mathbb{P}(\bar{A}(N; p_c, q_2)) \geq D_1 \mathbb{P}(\bar{A}(N; p_c, p_c))$$

and

$$\mathbb{P}(\bar{A}(N; q_1, p_c)) \leq D_2 \mathbb{P}(\bar{A}(N; p_c, p_c)).$$

Since the proofs of the above inequalities are similar, we only prove the first inequality. For that, we use a generalization of Russo's formula [5]. We take a small $\delta > 0$. The difference $\mathbb{P}(\bar{A}(N; p_c, p)) - \mathbb{P}(\bar{A}(N; p_c, p + \delta))$ can be written as the sum

$$\delta \sum_{e \in B(N), e \neq (1,0)} \mathbb{P}(\bar{A}(N; p_c, \cdot), \bar{A}_e(N; p, \cdot), D_e(N; p)) + O(\delta^2),$$

where $D_e(N; p)$ is the event that there exist three $p$-closed paths $P_1 - P_3$ in $B(N)^*$; the path $P_1$ connects an end of the edge $(1, 0)^*$ to an end of the edge $e^*$; the path $P_2$ connects the other end of the edge $(1, 0)^*$ to $R_N^*$ and the path $P_3$ connects the other end of the edge $e^*$ to $L_N^*$; or the path $P_2$ connects the other end of the edge $(1, 0)^*$ to $L_N^*$ and the path $P_3$ connects the other end of the edge $e^*$ to $R_N^*$. Letting $\delta$ tend to 0, we obtain

$$\frac{d}{dp} \mathbb{P}(\bar{A}(N; p_c, p)) = - \sum_e \mathbb{P}(\bar{A}(N; p_c, \cdot), \bar{A}_e(N; p, \cdot), D_e(N; p)). \tag{39}$$

We write the right-hand side of (39) as

$$
(40) \qquad -\sum_{j=1}^{[N/2]} \sum_{e\,:\,|e_x|=j} \mathbb{P}(\bar{A}(N;p_c,\cdot), \bar{A}_e(N;p,\cdot), D_e(N;p)),
$$

$$
(41) \qquad -\sum_{j=[N/2]+1}^{N} \sum_{e\,:\,|e_x|=j} \mathbb{P}(\bar{A}(N;p_c,\cdot), \bar{A}_e(N;p,\cdot), D_e(N;p)).
$$

By independence, the sum (40) is bounded from below by

$$
-\sum_{j=1}^{[N/2]} \sum_{e\,:\,|e_x|=j} \mathbb{P}(A([j/2];p_c,p))\mathbb{P}(A([3j/2],N;p_c,p))\mathbb{P}(A_e(e_x;[j/2];p,p)).
$$

We use (35), the bound $\sharp\{e : |e_x| = j\} \leq 16j$ and the fact that Lemma 6.3 is proved for $q_1 = q_2$ to bound the above sums from below by

$$
(42) \qquad \begin{aligned}
&-C_1 \sum_{j=1}^{[N/2]} j\mathbb{P}(A(j;p_c,p))\mathbb{P}(A(j,N;p_c,p))\mathbb{P}(A(j;p_c,p_c)) \\
&\geq -C_2\mathbb{P}(A(N;p_c,p)) \sum_{j=1}^{[N/2]} j\mathbb{P}(A(j;p_c,p_c)),
\end{aligned}
$$

where the inequality follows from (36). We estimate the sum in (42) using the relation

$$
(43) \qquad \sum_{j=1}^{N} j\mathbb{P}(A(j;p_c,p_c)) \asymp N^2 \mathbb{P}(A(N;p_c,p_c)).
$$

The relation (43) follows from (36) and the fact that $\mathbb{P}(A(j,N;p_c,p_c)) \geq C_3(j/N)^{2-C_4}$ for some positive $C_3$ and $C_4$ that do not depend on $j$ and $N$. This fact follows, for example, from [17], Theorem 24, where the 5-arms exponent is computed for site percolation on the triangular lattice. The same proof applies to bond percolation on the square lattice.

Similarly to the proof of [20], Lemma 6.2, the sum (41) can be bounded from below by

$$
-C_5 N^2 \mathbb{P}(A(N;p_c,p))\mathbb{P}(A(N;p_c,p_c)).
$$

This follows from a priori estimates of probabilities of two arms in a half-plane. We refer the reader to the proof of [20], Lemma 6.2, for more details. Again, although the proof of [20], Lemma 6.2, is given for site percolation on the triangular lattice, it also applies to bond percolation on the square lattice.

Putting together the bounds for the sums (40) and (41), and using (34), we obtain that the right-hand side of (39) is bounded from below by

$$-C_6 N^2 \mathbb{P}(\bar{A}(N; p_c, p)) \mathbb{P}(A(N; p_c, p_c)).$$

Therefore,

(44) $$\frac{d}{dp} \log \mathbb{P}(\bar{A}(N; p_c, p)) \geq -C_6 N^2 \mathbb{P}(A(N; p_c, p_c))$$

and

$$\mathbb{P}(\bar{A}(N; p_c, p)) \geq \mathbb{P}(\bar{A}(N; p_c, p_c)) e^{-C_6(p-p_c)N^2 \mathbb{P}(A(N; p_c, p_c))}$$
$$\geq \mathbb{P}(\bar{A}(N; p_c, p_c)) e^{-C_6(p_N-p_c)N^2 \mathbb{P}(A(N; p_c, p_c))}$$
$$\geq C_7 \mathbb{P}(\bar{A}(N; p_c, p_c)).$$

In the last inequality, we use (11). $\square$

DEFINITION 6.4. For any positive integers $n \leq m \leq 2m \leq N$ and edge $e \in \mathrm{Ann}(m, 2m)$, we define $C_e(n, N; m)$ as the event that:

- there exist two disjoint $p_c$-open paths $P_1$ and $P_2$ inside $\mathrm{Ann}(n, N) \setminus \{e\}$, the path $P_1$ connects $e_x$ or $e_y$ to $B(n)$ and the path $P_2$ connects the other end of $e$ to $\partial B(N)$; and
- there exists a $p_m$-closed path $P$ connecting $e_x^*$ and $e_y^*$ inside $\mathrm{Ann}(m, 2m)^* \setminus \{e^*\}$ so that $P \cup \{e^*\}$ is a circuit around the origin in $\mathrm{Ann}(m, 2m)^*$.

Note that if event $C_e(n, N; m) \cap \{\tau_e \in (p_c, p_m)\}$ occurs, then there is no $p_c$-open crossing of $\mathrm{Ann}(n, N)$ and no $p_m$-closed circuit in $\mathrm{Ann}(m, 2m)^*$ (see Figure 3).

DEFINITION 6.5. Let $n$, $m$ and $N$ be positive integers such that $2n \leq m$ and $3m \leq N$. Let $x = ([m/2], [3m/2])$. For $e \in B(x, [m/2])$, we define $\widetilde{C}_e(n, N; m)$ as the event that:

- the event $\widetilde{A}_e(x; m; p_c, p_m)$ occurs;
- there are two disjoint $p_c$-open paths $P_5$ and $P_6$ such that $P_5$ connects $U_{[m/2]}(x)$ to the boundary of $B(N)$ inside $\mathrm{Ann}(2m-1, N)$ and $P_6$ connects $D_{[m/2]}(x)$ to the boundary of $B(n)$ inside $\mathrm{Ann}(n, m)$. Moreover, $P_5$ and $P_6$ satisfy $P_5 \cap \mathrm{Ann}(x; [m/2], m) \subset U_{[m/2], m}(x)$ and $P_6 \cap \mathrm{Ann}(x; [m/2], m) \subset D_{[m/2], m}(x)$;
- there exists a $p_m$-closed path $P$ inside $\mathrm{Ann}(m, 2m-1)^* \setminus B(x, [m/2])^*$ such that $P$ connects $L_{[m/2]}(x)^*$ to $R_{[m/2]}(x)^*$ and $P \cap \mathrm{Ann}(x; [m/2], m)^* \subset L_{[m/2], m}(x)^* \cup R_{[m/2], m}(x)^*$.

The event $\widetilde{C}_e(n, N; m)$ is illustrated in Figure 4.

The event $\widetilde{C}_e(n, N; m)$ obviously implies the event $C_e(n, N; m)$. The reason we introduce the event $\widetilde{C}_e(n, N; m)$ is that

$$(45) \qquad \mathbb{P}(\widetilde{C}_e(n, N; m)) \asymp \mathbb{P}(\widetilde{A}_e(x; m; p_c, p_m))\mathbb{P}_{cr}(B(n) \leftrightarrow \partial B(N)),$$

where the constants do not depend on $e$, $m$, $n$ and $N$. This observation follows from Lemma 6.1, the RSW theorem, and (35) and (36) applied to $q_1 = q_2 = p_c$.

COROLLARY 6.2. *For any positive integers* $n$, $m$ *and* $N$ *such that* $2n \leq m$ *and* $3m \leq N$,

$$(46) \qquad \begin{aligned} &\mathbb{P}(\exists e \in \mathrm{Ann}(m, 2m): \tau_e \in (p_c, p_m), C_e(n, N; m)) \\ &\qquad \geq C_8 \mathbb{P}_{cr}(B(n) \leftrightarrow \partial B(N)), \end{aligned}$$

*where* $C_8$ *does not depend on* $n$, $N$ *and* $m$.

PROOF. Note that the events

$$\{\tau_e \in (p_c, p_m), C_e(n, N; m)\}_{e \in \mathrm{Ann}(m, 2m)}$$

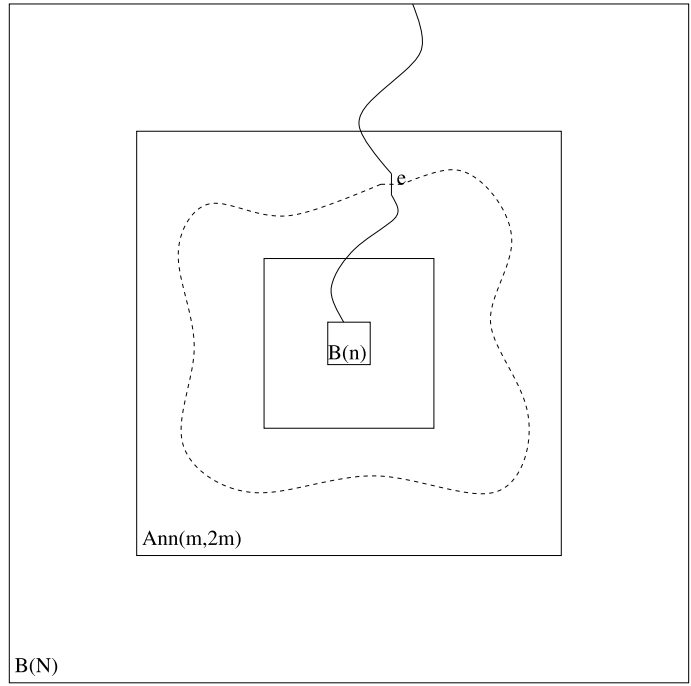

FIG. 3. *Event* $C_e(n, N; m)$. *The solid curves represent* $p_c$-*open paths and the dotted curves represent* $p_m$-*closed paths. The edge* $e$ *does not have to be* $p_c$-*open or* $p_m$-*closed.*

are disjoint. Therefore,

$$\mathbb{P}(\exists e \in \mathrm{Ann}(m, 2m) : \tau_e \in (p_c, p_m), C_e(n, N; m))$$

$$= \sum_{e \in \mathrm{Ann}(m, 2m)} \mathbb{P}(\tau_e \in (p_c, p_m), C_e(n, N; m))$$

$$\geq (p_m - p_c) \sum_{e \in B(x, [m/2])} \mathbb{P}(\widetilde{C}_e(n, N; m))$$

$$\geq C_9(p_m - p_c) \sum_{e \in B(x, [m/2])} \mathbb{P}(\widetilde{A}_e(x; m; p_c, p_m)) \mathbb{P}_{cr}(B(n) \leftrightarrow \partial B(N))$$

$$\geq C_{10}(p_m - p_c) m^2 \mathbb{P}(A(m; p_c, p_c)) \mathbb{P}_{cr}(B(n) \leftrightarrow \partial B(N))$$

$$\geq C_{11} \mathbb{P}_{cr}(B(n) \leftrightarrow \partial B(N)).$$

The second inequality follows from (45). In the third inequality, we use (37) and Lemma 6.3. In the last inequality, we use (11).  □

PROOF OF THEOREM 1.5. LOWER BOUND. We give the proof for $k = 2$. The case $k = 1$ was considered in [18] and the proof for $k \geq 3$ is similar to the one for $k = 2$. Note that the event $\{\widehat{R}_2 > n\}$ is implied by the event that there exists an edge $e \in B(n)$ and $p > p_c$ such that:

- $\tau_e \in (p_c, p)$;
- there exist two $p_c$-open paths $P_1$ and $P_2$ in $B(n)$, the path $P_1$ connects an end of $e$ to the origin and the path $P_2$ connects the other end of $e$ to the boundary of $B(n)$;
- there exists a $p$-closed path $P$ in $B(n)^*$ connecting $e_x^*$ to $e_y^*$ so that $P \cup \{e^*\}$ is a circuit around the origin.

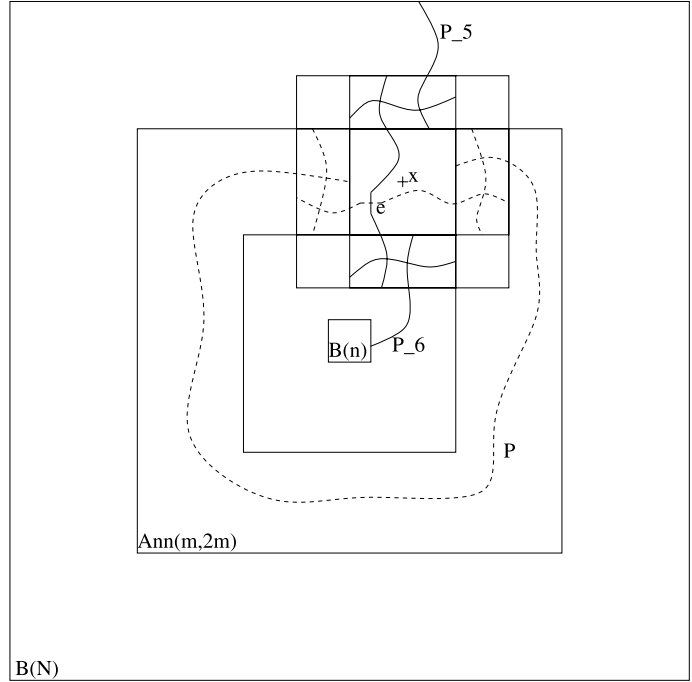

FIG. 4. Event $\widetilde{C}_e(n, N; m)$. The solid curves represent $p_c$-open paths and the dotted curves represent $p_m$-closed paths. The edge $e$ does not have to be $p_c$-open or $p_m$-closed.

There could be at most one edge $e \in B(n)$ which satisfies the above three conditions. Therefore,

$$\mathbb{P}(\widehat{R}_2 > n) \geq \sum_{k=0}^{[\log n]-1} \mathbb{P}(\exists e \in \mathrm{Ann}([n/2^{k+1}], [n/2^k]) : \tau_e \in (p_c, p_{n/2^{k+1}}),$$

$$C_e(1, n; [n/2^{k+1}]))$$

$$\geq C_{12} \sum_{k=0}^{[\log n]-1} \mathbb{P}_{cr}(0 \leftrightarrow \partial B(n))$$

$$= C_{12}[\log n] \mathbb{P}_{cr}(0 \leftrightarrow \partial B(n)).$$

The last inequality follows from (46). $\quad\square$

**7. Proof of Theorem 1.8.** Let $G = (\mathcal{G}, \mathcal{E})$ be an infinite connected subgraph of $(\mathbb{Z}^2, \mathbb{E}^2)$ which contains the origin. We call an edge $e \in \mathcal{E}$ a *disconnecting edge* for $G$ if the graph $(\mathcal{G}, \mathcal{E} \setminus \{e\})$ has a finite component and if the origin belongs to this finite component. Note that each outlet of the invasion is a disconnecting edge for the IPC.

Let $D_{m,n}$ be the event that the IIC does not contain a disconnecting edge in the annulus $\mathrm{Ann}(m, n)$ and let $\mathcal{D}_{m,n}$ be the event that the IPC does not contain a disconnecting edge in the annulus $\mathrm{Ann}(m, n)$. We prove the following theorem:

THEOREM 7.1. *There exists a sequence $(n_k)$ such that*

$$(47) \qquad \mathbb{P}\left( \sum_k \mathbb{I}(\mathcal{D}_{n_k, n_{k+1}}) < \infty \right) = 1$$

*and*

$$(48) \qquad \nu\left( \sum_k \mathbb{I}(D_{n_k, n_{k+1}}) = \infty \right) = 1.$$

Theorem 1.8 immediately follows from Theorem 7.1. Indeed, Theorem 7.1 implies that the IIC is supported on clusters for which infinitely many of the events $D_{n_k,n_{k+1}}$ occur and the IPC is supported on clusters for which only finitely many of the events $\mathcal{D}_{n_k,n_{k+1}}$ occur. Roughly speaking, this says that the distance between consecutive disconnecting edges (ordered by distance from the origin) can be much larger in the IIC than in the IPC. The proof of Theorem 7.1 is based on the following result (see Section 1.2 for the definitions).

THEOREM 7.2. *There exist $C_1, C_2$ such that for all $1 \leq m < n$,*

$$\mathbb{P}(\mathcal{D}_{m,n}) \leq C_1 \mathbb{P}_{cr}(A^2_{m,n}) \tag{49}$$

*and*

$$\nu(D_{m,n}) \geq C_2 \frac{\mathbb{P}_{cr}(A^2_{m,n})}{\mathbb{P}_{cr}(A^1_{m,n})}. \tag{50}$$

LEMMA 7.1 ([17], Theorem 27). *For all positive integers $m < n$ and for all $p \in [p_c, p_n]$,*

$$\mathbb{P}_p(A^2_{m,n}) \asymp \mathbb{P}_{cr}(A^2_{m,n}),$$

*where the constants do not depend on $m$, $n$ and $p$.*

Although Theorem 27 in [17] is stated for site percolation on the triangular lattice, the proof for bond percolation on the square lattice is the same.

LEMMA 7.2. *There exists $C_3$ such that for all $m_1 < m_2 < n$, we have*

$$\frac{\mathbb{P}_{cr}(A^2_{m_1,n})}{\mathbb{P}_{cr}(A^2_{m_1,m_2})} \geq C_3 \frac{m_2}{n}.$$

PROOF. This follows from a priori estimates of probabilities of two arms in a half-plane (see [17], Theorem 24). □

PROOF OF THEOREM 7.2. We first prove (49). Note that if the invasion percolation cluster contains a circuit, then there is a pond that entirely contains this circuit. Therefore, the event $\mathcal{D}_{m,n}$ can only occur if there exists an invasion pond which contains two disjoint crossings $P_1$ and $P_2$ of the annulus $\mathrm{Ann}(m,n)$ (see Figure 5). Therefore, there exists $p'$ such that $P_1$ and $P_2$ are $p'$-open and there exists a circuit around the origin which is $p'$-closed and which has diameter at least $n$.

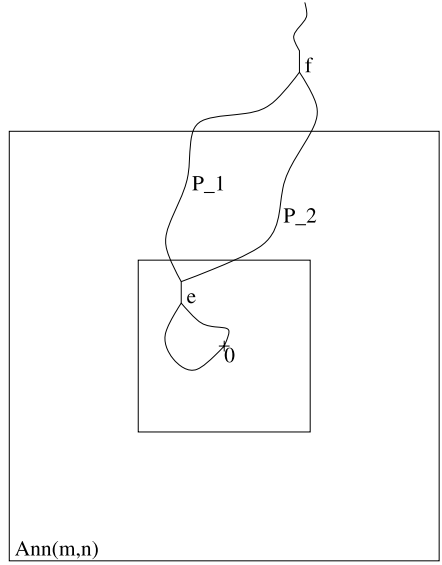

Fig. 5. *Event $\mathcal{D}_{m,n}$. The edges $e$ and $f$ are disconnecting. The paths $P_1$ and $P_2$ create a circuit, which implies that there is a pond that entirely contains both paths.*

Recall the definition of $(p_n(j))$ from (5). Later, we take $C_*$ in (5) to be sufficiently large. We decompose the event $\mathcal{D}_{m,n}$ according to the value of $p'$:

$$(51) \qquad \mathbb{P}(\mathcal{D}_{m,n}) = \sum_{j=1}^{\log^* n} \mathbb{P}(\mathcal{D}_{m,n}; p' \in [p_n(j), p_n(j-1))).$$

Note that the event $\{\mathcal{D}_{m,n}; p' \in [p_n(j), p_n(j-1))\}$ implies the event $A^2_{m,n,p_n(j-1)} \cap B_{n,p_n(j)}$ (see Section 1.2 for the definition of these events). It follows from (6) and (10) that there exist constants $C_4$ and $C_5$ such that the probability $\mathbb{P}(B_{n,p_n(j)})$ is bounded from above by $C_4(\log^{(j-1)} n)^{-C_5}$. The constant $C_5$ can be made arbitrarily large by making $C_*$ large enough. We use Lemmas 7.1 and 7.2 to bound the probability $\mathbb{P}(A^2_{m,n,p_n(j-1)}) \leq C_6(\log^{(j-1)} n)\mathbb{P}_{cr}(A^2_{m,n})$. We use the FKG inequality and the above estimates for the events $A^2_{m,n,p_n(j-1)}$ and $B_{n,p_n(j)}$ to get

$$\mathbb{P}(\mathcal{D}_{m,n}) \leq C_4 C_6 \mathbb{P}_{cr}(A^2_{m,n}) \sum_{j=1}^{\log^* n} (\log^{(j-1)} n)^{1-C_5}$$

$$\leq C_7 \mathbb{P}_{cr}(A^2_{m,n}).$$

The last inequality follows from [6], (2.26), if we take $C_*$ such that $C_5 > 1$.

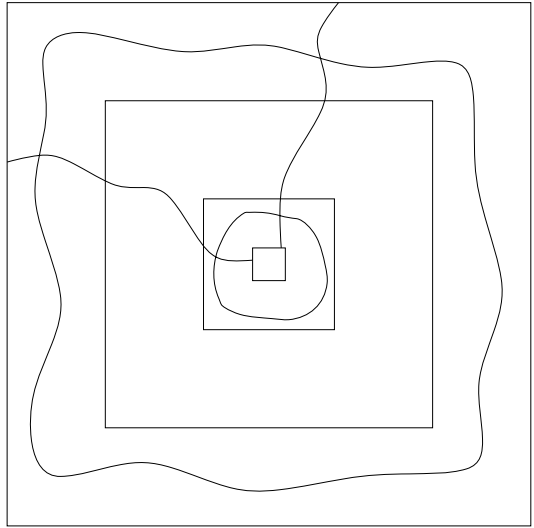

Fig. 6. *Event* $C_{m,n}$. *The inner circuit is in the annulus* Ann($[m/2], m$) *and the outer circuit is in the annulus* Ann($n, 2n$).

We now prove (50). Let $C_{m,n}$ be the event $A_{[m/2],m} \cap A_{[m/2],2n}^2 \cap A_{n,2n}$ (see Figure 6). Note that $\nu(C_{m,n} \setminus D_{m,n}) = 0$. It is therefore sufficient to prove (50) for $C_{m,n}$.

For positive integers $m < n < N$ (later, we consider the limit as $N$ tends to infinity), we use the FKG inequality to get

$$\mathbb{P}_{cr}(C_{m,n} \cap A_{0,N}^1)$$
$$\geq \mathbb{P}_{cr}(C_{m,n} \cap A_{0,m}^1 \cap A_{n,N}^1)$$
$$\geq C_8 \mathbb{P}_{cr}(A_{[m/2],2n}^2) \mathbb{P}_{cr}(A_{0,m}^1) \mathbb{P}_{cr}(A_{n,N}^1)$$
$$\geq C_8 \frac{\mathbb{P}_{cr}(A_{[m/2],2n}^2) \mathbb{P}_{cr}(A_{0,N}^1)}{\mathbb{P}_{cr}(A_{m,n}^1)}$$

for some $C_8 > 0$. Standard RSW arguments give a constant $C_9$ such that for all $1 \leq m < n$,

$$\mathbb{P}_{cr}(A_{[m/2],2n}^2) \geq C_9 \mathbb{P}_{cr}(A_{m,n}^2).$$

Therefore,

$$\nu(D_{m,n}) \geq C_8 C_9 \frac{\mathbb{P}_{cr}(A_{m,n}^2)}{\mathbb{P}_{cr}(A_{m,n}^1)}. \qquad \qquad \square$$

PROPOSITION 7.1. *There exists a sequence $(n_k)$ such that $n_{k+1} > 4n_k$,*

$$\sum_k \mathbb{P}_{cr}(A^2_{n_k, n_{k+1}}) < \infty \tag{52}$$

*and*

$$\sum_k \frac{\mathbb{P}_{cr}(A^2_{n_{2k}, n_{2k+1}})}{\mathbb{P}_{cr}(A^1_{n_{2k}, n_{2k+1}})} = \infty. \tag{53}$$

Proposition 7.1 follows from Lemma 7.2 and the fact that $\mathbb{P}_{cr}(A^1_{m,n}) \leq c(m/n)^\delta$ for some positive $c$ and $\delta$. Indeed, we obtain $\mathbb{P}_{cr}(A^1_{m,n}) \leq c(C_3\mathbb{P}_{cr}(A^2_{m,n}))^\delta$. We now take, for example, the sequence $n_k = \min\{n > 4n_{k-1} : \mathbb{P}_{cr}(A^2_{n_{k-1},n}) \leq (1/k)^{1+\delta}\}$.

PROOF OF THEOREM 7.1. We take a sequence from Proposition 7.1. Equality (47) follows from the Borel–Cantelli lemma. To prove (48), we use Borel's lemma [13]:

LEMMA 7.3. *Consider a probability space $(\Omega, \mathcal{F}, \mathbb{P})$ and a sequence of events $\Gamma_n \in \mathcal{F}$. Let $\limsup_n \Gamma_n = \bigcap_n \bigcup_{k \geq n} \Gamma_k$ be the event that infinitely many of the $\Gamma_n$'s occur. Let $a_n = \mathbb{I}(\Gamma_n)$ be the indicator of event $\Gamma_n$. If there exists a sequence $b_n$ such that $\sum_n b_n = \infty$ and for any $\alpha_i \in \{0, 1\}$, $i = 1, \ldots, n-1$,*

$$\mathbb{P}(\Gamma_n | a_1 = \alpha_1, \ldots, a_{n-1} = \alpha_{n-1}) \geq b_n > 0,$$

*then*

$$\mathbb{P}\left(\limsup_n \Gamma_n\right) = 1.$$

Note that it is sufficient to prove (48) for the events $C_{n_k, n_{k+1}}$ (see the proof of Theorem 7.2 for the definition). We apply Lemma 7.3 to the probability measure $\nu$ and to the events $C_{n_{2k}, n_{2k+1}}$. Let $d_k = \mathbb{I}(C_{n_{2k}, n_{2k+1}})$. A slight extension of the proof of (50) gives, for any $\alpha_i \in \{0, 1\}$, $i = 1, \ldots, k-1$,

$$\nu(C_{n_{2k}, n_{2k+1}} | d_1 = \alpha_1, \ldots, d_{k-1} = \alpha_{k-1}) \geq C_2 \frac{\mathbb{P}_{cr}(A^2_{n_{2k}, n_{2k+1}})}{\mathbb{P}_{cr}(A^1_{n_{2k}, n_{2k+1}})} =: b_k, \tag{54}$$

where $C_2$ is the constant from (50). Indeed, let $\mathcal{W}$ be the set of configurations of edges in $B(2n_{2k-1})$ such that $d_1 = \alpha_1, \ldots, d_{k-1} = \alpha_{k-1}$. For any $\omega \in \mathcal{W}$ and large enough $N$,

$$\mathbb{P}_{cr}(C_{n_{2k}, n_{2k+1}} \cap A^1_{0,N} | \omega)$$
$$\geq \mathbb{P}_{cr}(C_{n_{2k}, n_{2k+1}} \cap A^1_{0, n_{2k}} \cap A^1_{n_{2k+1}, N} | \omega)$$

$$\geq C_8 \mathbb{P}_{cr}(A^2_{[n_{2k}/2],2n_{2k+1}}|\omega)\mathbb{P}_{cr}(A^1_{0,n_{2k}}|\omega)\mathbb{P}_{cr}(A^1_{n_{2k+1},N}|\omega)$$

$$= C_8 \mathbb{P}_{cr}(A^2_{[n_{2k}/2],2n_{2k+1}})\mathbb{P}_{cr}(A^1_{0,n_{2k}}|\omega)\mathbb{P}_{cr}(A^1_{n_{2k+1},N})$$

$$\geq C_8 \frac{\mathbb{P}_{cr}(A^2_{[n_{2k}/2],2n_{2k+1}})\mathbb{P}_{cr}(A^1_{0,N}|\omega)}{\mathbb{P}_{cr}(A^1_{n_{2k},n_{2k+1}})},$$

which implies (54). In the second line, we used the FKG inequality and independence. The equality follows from independence. From the choice of $(n_k)$, it follows that $\sum_k b_k = \infty$. Therefore, equality (48) follows from Lemma 7.3. $\square$

**Acknowledgments.** We would like to thank Rob van den Berg and Antal Járai for suggesting these problems. We thank Rob van den Berg and Federico Camia for enjoyable discussions. Finally, we are indebted to Rob van den Berg for careful readings of the manuscript and for many useful comments.

M. DAMRON                                A. SAPOZHNIKOV
COURANT INSTITUTE                        EURANDOM
 OF MATHEMATICAL SCIENCES                P.O. BOX 513
251 MERCER STREET                        5600 MB EINDHOVEN
NEW YORK, NEW YORK 10012                 THE NETHERLANDS
USA                                      E-MAIL: sapozhnikov@eurandom.tue.nl
E-MAIL: damron@cims.nyu.edu

              B. VÁGVÖLGYI
              VRIJE UNIVERSITEIT AMSTERDAM
              DE BOELELAAN 1081
              1081 HV AMSTERDAM
              THE NETHERLANDS
              E-MAIL: bvagvol@few.vu.nl